\documentclass[sn-mathphys-num]{sn-jnl}

\usepackage{graphicx}%
\usepackage{multirow}%
\usepackage{amsmath,amssymb,amsfonts}%
\usepackage{amsthm}%
\usepackage{mathrsfs}%
\usepackage[title]{appendix}%
\usepackage{xcolor}%
\usepackage{textcomp}%
\usepackage{manyfoot}%
\usepackage{booktabs}%
\usepackage{algorithm}%
\usepackage{algorithmicx}%
\usepackage{algpseudocode}%
\usepackage{listings}%
\usepackage{color,xcolor}
\usepackage{amsmath,amssymb,bm}
\usepackage{graphicx,graphics}
\definecolor{refgreen}{rgb}{0,0.5,0}
\usepackage{dsfont}
\usepackage{yfonts}

\usepackage{tikz}
\usetikzlibrary{arrows}
\usepackage{adjustbox}

\usepackage{leftidx}

\usepackage{amsthm}
\newtheorem{thm}{Theorem}[section]
\newtheorem{lem}[thm]{Lemma}
\newtheorem{prop}[thm]{Proposition}
\newtheorem{rem}[thm]{Remark}

\newcommand\bfd{{\mathbf d}}

\newcommand\bfs{{\mathbf s}}
\newcommand\bfu{{\mathbf u}}

\newcommand\bfpsi{{\boldsymbol \psi}}

\newcommand\calT{{\mathcal T}}
\newcommand\calE{{\mathcal E}}

\newcommand\calI{{\mathcal I}}
\def            \d          {\hspace{1.5pt}\mathrm{d}}

\newcommand{\R}{\mathbb{R}}
\DeclareMathOperator{\Span}{span}
\DeclareMathOperator{\Real}{Re}

\newcommand{\Ga}{\varGamma}
\newcommand{\Om}{\varOmega}
\renewcommand{\nu}{\textnormal{n}}
\newcommand{\matdev}{\partial^{\bullet}}

\newcommand{\h}{\mathds{h}}
\newcommand{\m}{\mathbf{m}}

\newcommand{\mh}{m_{\Ga_h}}
\newcommand{\ah}{a_{\Ga_h}}

\newcounter{scheme}[section]

\begin{document}
	
	\title[Error estimates for full discretization by an almost mass conservation technique for Cahn--Hilliard systems with dynamic boundary conditions]{Error estimates for full discretization by an almost mass conservation technique for Cahn--Hilliard systems with dynamic boundary conditions}
	
	
	\author*[1]{\fnm{Nils} \sur{Bullerjahn}}\email{bullerja@uni-paderborn.de}
	
	
	
	\affil*[1]{\orgdiv{Institute of Mathematics}, \orgname{Paderborn University}, \orgaddress{\street{Warburgerstr.~100}, \city{Paderborn}, \postcode{33098},  \country{Germany}}}
	
	
	
	
	\abstract{A proof of optimal-order error estimates is given for the full discretization of the bulk--surface Cahn--Hilliard system with dynamic boundary conditions in a smooth domain. The numerical method combines a linear bulk--surface finite element discretization in space and linearly implicit backward difference formulae of order one to five in time. The error estimates are obtained by a consistency and stability analysis, based on an energy estimate and the novel approach of exploiting the almost mass conservation of the error equations to derive a Poincar\'e-type inequality. We demonstrate how this approach can be generalized to other almost mass conserving problems. To this end we prove optimal-order fully discrete error estimates for the Cahn--Hilliard equation on evolving surfaces. We illustrate and complement our findings by numerical experiments.
	}

	
	
	
	\maketitle
	
	\section{Introduction}
	In this paper we prove optimal-order error estimates for a full discretization of the bulk--surface Cahn--Hilliard system with dynamic boundary conditions---a two parameter family of dynamic boundary conditions generalizing the most important dynamic boundary conditions of Cahn--Hilliard-type. We present, as the main focus of this paper, a novel technique for achieving optimal-order error estimates for the full discretization of the bulk--surface Cahn--Hilliard system with dynamic boundary conditions using linear bulk--surface finite elements in space and linearly implicit backward difference method of order $1$ to $5$ in time. 
	
	The new idea for the stability analysis uses mass conservation properties of the equation to conclude an \emph{almost mass conservation} of the error, and combines this with a \emph{Poincar\'e-type estimate} and an energy estimate. We prove optimal-order error estimates, between the solution of the numerical scheme and a sufficiently regular exact solution, for the $L^2$- and $H^1$-norm of the phase field in bulk and surface, respectively. 
	
	This novel technique can be transferred to a wider class of equations, satisfying an energy estimate and mass conservation properties. In order to demonstrate this, we prove optimal-order fully discrete error estimates for the Cahn-Hilliard equation on evolving surface, discretized by evolving surface finite elements in space and linearly implicit backward difference method in time.

	\smallskip
	The Cahn--Hilliard equation was developed as a diffuse interface model for the phase separation process of binary mixtures in a given domain by \cite{CahnHilliard1958}. Originally homogeneous Neumann boundary conditions were chosen to show well-posedness of the problem, however this choice leads to some limitations of the model, for example the contact angle between the free interface of the two components and the boundary of the domain is fixed at $\pi/2$, see, e.g. \cite{Jacqmin2000}. In order to correctly model the behavior of the system at the boundary of the domain in accordance with different physical conditions, various dynamic boundary conditions were developed, see e.g. \cite{FischerMaass1997} for Allen--Cahn-type dynamic boundary conditions, \cite{GoldsteinMiranvilleSchimperna2011} and \cite{LiuWu2019} for Cahn--Hilliard-type dynamic boundary conditions. The extension of the last model by a Robin boundary approximation of an affine linear transmission condition, to the Cahn--Hilliard system with transmission rate dependent dynamic boundary conditions was developed in \cite{KnopfLam2020} and a similar condition was introduced in \cite{KnopfLamLiuMetzger2021} for the chemical potential, which leads to reaction rate dependent dynamic boundary conditions, interpolating between the previously mentioned Cahn--Hilliard type dynamic boundary conditions. In this paper we consider the combination of these Robin boundary extensions in full generality to the bulk--surface Cahn--Hilliard system with dynamic boundary conditions, recently introduced in \cite{KnopfStange2024}, which contains all the previously mentioned models \cite{GoldsteinMiranvilleSchimperna2011,LiuWu2019,KnopfLam2020,KnopfLamLiuMetzger2021} as special cases, cf. Figure~\ref{fig:graphDynBC}. A nice review paper collecting results on the Cahn--Hilliard equation with dynamic boundary conditions up to 2022 is \cite{Wu2022}.
	
	\smallskip 
	There are several works that deal with the Cahn--Hilliard equation \emph{with} dynamic boundary conditions, analyzing \emph{semi-discrete} convergence rates, fully discrete convergence without rates or qualitative properties. In \cite{CherfilsPetcuPierre2010} and \cite{CherfilsPetcu2014}, the authors have considered a backward Euler time discretization of Allen--Cahn-type, respectively, Cahn--Hilliard dynamic boundary conditions and derived optimal-order semi-discrete error estimates in space on a 2$d$ or 3$d$ slab. In \cite{HarderKovacs2022} there were shown optimal-order semi-discrete error estimates for the model in \cite{GoldsteinMiranvilleSchimperna2011} and a bulk--surface finite element approximation. The weak convergence up to a subsequence of the Cahn--Hilliard equation with reaction rate dependent dynamic boundary conditions for a bulk--surface finite element/backward Euler scheme was proved in \cite{KnopfLamLiuMetzger2021,Metzger2021}, and in \cite{Metzger2023} for an SAV scheme. In \cite{BaoZhang2021} the authors considered a backward Euler time discretization of reaction rate dependent dynamic boundary conditions and derived first order semi-discrete error estimates in time. Numerical experiments for the Cahn--Hilliard system with transmission rate dependent dynamic boundary conditions are presented in \cite{KnopfLam2020}, with a discretization by the backward Euler method in time and linear finite elements in space. For the dynamic boundary conditions proposed in \cite{LiuWu2019} and a semi-discrete in time stabilized linearly implicit BDF2 scheme, an energy dissipation and semi-discrete error estimates in time are shown in \cite{MengBaoZhang2023}, and in \cite{LiuShenZheng2024} this scheme is used for various dynamic boundary conditions in a multiple scalar auxiliary variables (MSAV) approach. In \cite{AltmannZimmer2023} the authors derive a dissipation-preserving scheme by formulating a partial differential--algebraic equation (PDAE) for Allen--Cahn-type, and Cahn--Hilliard-type dynamic boundary conditions, and a semi-discrete in time backward Euler approach with convex--concave splitting. Additional structure-preserving schemes for Cahn--Hilliard-type dynamic boundary conditions can be found in \cite{OkumuraFukao2024}. 
	
	To our knowledge the first and so far only fully discrete error estimates are provided in the recent work \cite{BullerjahnKovacs2024}, for the special case of the dynamic boundary conditions developed in \cite{GoldsteinMiranvilleSchimperna2011}, extending \cite{HarderKovacs2022} to optimal-order fully discrete error estimates for a $q$-step backward difference formulae (BDF method) of order $1$ to $5$. 
	
	The contribution of the present work are optimal-order fully discrete error estimates for almost the entire family of Cahn--Hilliard-type dynamic boundary conditions developed in \cite{KnopfStange2024}, including the models in \cite{GoldsteinMiranvilleSchimperna2011,KnopfLam2020,KnopfLamLiuMetzger2021}, cf. Figure~\ref{fig:graphDynBC}. 
	
	\smallskip
	The fundamental technique of analyzing an isoparametric bulk--surface finite element method for an elliptic partial differential equation in a curved domain was introduced in \cite{ElliottRanner2013}, and many general results are proved and summarized in \cite{ElliottRanner2021}, even for a possibly evolving domain. General linear and semi-linear parabolic equations of second order with dynamic boundary conditions were analyzed in \cite{KovacsLubich2017}, providing optimal-order error estimates for a large class of second-order problems in an abstract setting. 
	The spatial discretization and its error estimates for our scheme are based on the ideas developed in these papers, and the consistency analysis of the spatial defects will heavily depend on the general geometric approximation error estimates developed therein. 
	Some early results for conforming finite element discretizations of parabolic problems with dynamic boundary conditions on polyhedral domains were shown in \cite{Fairweather}.
	
	\smallskip 
	In order to adapt the energy estimates for the $q$-step BDF method, we use a combination of the results from the $G$-stability theory of \cite{Dahlquist1978} and the multiplier technique of \cite{NevanlinnaOdeh1981}. The approach of testing with the errors was first used for PDEs in \cite{LubichMansourVenkataraman2013}, for linear parabolic differential equations on evolving surfaces, and later on for abstract quasi-linear parabolic problems in \cite{AkrivisLubich2015}. Fully discrete energy estimates---testing with the discrete time derivatives of the errors---were developed for mean curvature flow in \cite{KovacsLiLubich2019} and recently used in \cite{BullerjahnKovacs2024} for the Cahn--Hilliard equation with Cahn--Hilliard-type dynamic boundary conditions. The present work requires similar arguments for testing with the errors as well as their discrete time derivatives.	
	
	\smallskip
	The key novelty of this paper is the introduction of a discrete Poincar\'e-type inequality, based on the \emph{almost mass conservation} of the error equations of the numerical scheme. This offers a new tool in analyzing this and similar equations by BDF fully discrete energy estimates and eliminates the dependence on the anti-symmetric structure of the arguments in \cite{BullerjahnKovacs2024}. It is therefore to be expected that these techniques can be transferred to other similar equations, such as Allen--Cahn and Cahn--Hilliard equations in stationary domains with classical boundary conditions, see \cite{AllenCahn1979}, respectively \cite{CahnHilliard1958}, or the same PDEs on evolving surfaces, see \cite{OlshanskiiXuYushutin2021}, respectively \cite{CaetanoElliott2021}. We demonstrate this by applying the strategy of proof to the Cahn--Hilliard equation on evolving surfaces, proving optimal-order error estimates for the fully discrete numerical scheme. Comparable results were obtained in \cite{ElliottSales2024}, \cite{ElliottSales2024B}, where the authors derived optimal-order error estimates for an evolving surface finite element scheme and a fully implicit backward Euler time-discretization, for a smooth potential with polynomial growth, respectively logarithmic potential. In contrast to these results we are able to show optimal-order error estimates for backward difference formula up to order $5$, with a potential only satisfying local Lipschitz continuity for up to three of its derivatives, and even super-convergence in the $H^1$-norm when comparing with the Ritz map of the exact solution.
	
	\smallskip
	The paper is structured as follows:
	In Section~\ref{section:Cahn--Hilliard intro} we introduce the bulk--surface Cahn--Hilliard system with dynamic boundary conditions, which will be the basis of our discretization and in Section~\ref{section:Notation} we formulate some notation, preliminaries and assumptions.
	In Section~\ref{section:Discretization CH} we introduce the full discretization in space by bulk--surface finite elements, and in time by the $q$-step backward difference method and in Section~\ref{section:MainResult} we state and discuss our main convergence result, while in Section~\ref{section:BS - Poincare-wirtinger} we state the bulk--surface Poincar\'e--Wirtinger inequality. In Section~\ref{section:Stability} we establish the stability estimates, Section~\ref{section:Consistency} is devoted to the consistency analysis, and Section~\ref{section:MainProof} will combine these results to prove our main result.
	Finally in Section~\ref{section:Generalization} we demonstrate how these strategies can be adapted to other mass conserving problems, by proving optimal-order fully discrete error estimates for the Cahn--Hilliard equation on evolving surfaces. In Section~\ref{section:NumericalExperiments} we present some numerical experiments to illustrate and complement our theoretical results.
	
	\section{Bulk--surface Cahn--Hilliard system with dynamic boundary conditions \label{section:Cahn--Hilliard intro}}
	
\begin{figure}[hb] 
		\begin{adjustbox}{width=\textwidth}
			\begin{tikzpicture}
				\draw[-stealth, very thick, dashed] (0,0)  -- (5.4,0) node[right] {$K$}; 
				\node at (0,-0.3) {$0$};
				\node at (5,-0.3) {$\infty$};
				\draw[-stealth, very thick] (0,0)  -- (0,5.4) node[above] {$L$}; 
				\node at (-0.3,0) {$0$};
				\node at (-0.3,5) {$\infty$};
				
				\draw[dashed, very thick] (5,0) -- (5,5);
				\draw[dashed, very thick] (0,5) -- (5,5);
				
				\node[draw,align=center] at (-1,5.6) {\footnotesize LW-model \\ \footnotesize \cite{LiuWu2019}};
				
				\node[draw,align=center] at (2.5,5.6) {\footnotesize transmission rate dependent \\ \footnotesize dyn. bound. cond. \cite{KnopfLam2020}};
				
				\node[draw,align=center] at (-1.1,-0.6) {\footnotesize GMS-model \\ \footnotesize \cite{GoldsteinMiranvilleSchimperna2011}};
				
				\node[draw,rotate=90, align=center] at (-0.6,2.5) {\footnotesize reaction rate dependent \\ \footnotesize dyn. bound. cond. \cite{KnopfLamLiuMetzger2021}};
				
				\node[draw,align=center] at (6.6,5.6) {\footnotesize Neumann boundary \\ \footnotesize conditions \cite{CahnHilliard1958}};

				\fill[gray!20] (0.01,0.01) rectangle (4.99,4.99);
				
				\draw[gray!50, very thick] (0,5) circle (0.1);
				\filldraw[color=gray!50, fill=gray!50, very thick] (0,0) circle (0.1);
				\filldraw[color=gray!50, fill=gray!50, very thick] (5,5) circle (0.1);
				
				\node[draw,align=center] at (2.5,2.5) {\footnotesize bulk--surface system with \\ \footnotesize dyn. bound. cond. \cite{KnopfStange2024}};
				
				\draw[gray!50, very thick] (0.1,5) -- (5,5);
				\draw[gray!50, very thick] (0,0.1) -- (0,4.9);
				
			\end{tikzpicture}
		\end{adjustbox}
		\caption{Sketch of the different types of dynamic boundary conditions corresponding to the parameter choices $K,L \in [0,\infty]$. The main result of this paper, Theorem~\ref{thm:OptimalOrderErrorEst} shows optimal-order error estimates for the range of parameters marked in gray (the domain $[0,\infty]^2$ except circles and dashed lines).}\label{fig:graphDynBC}
	\end{figure}
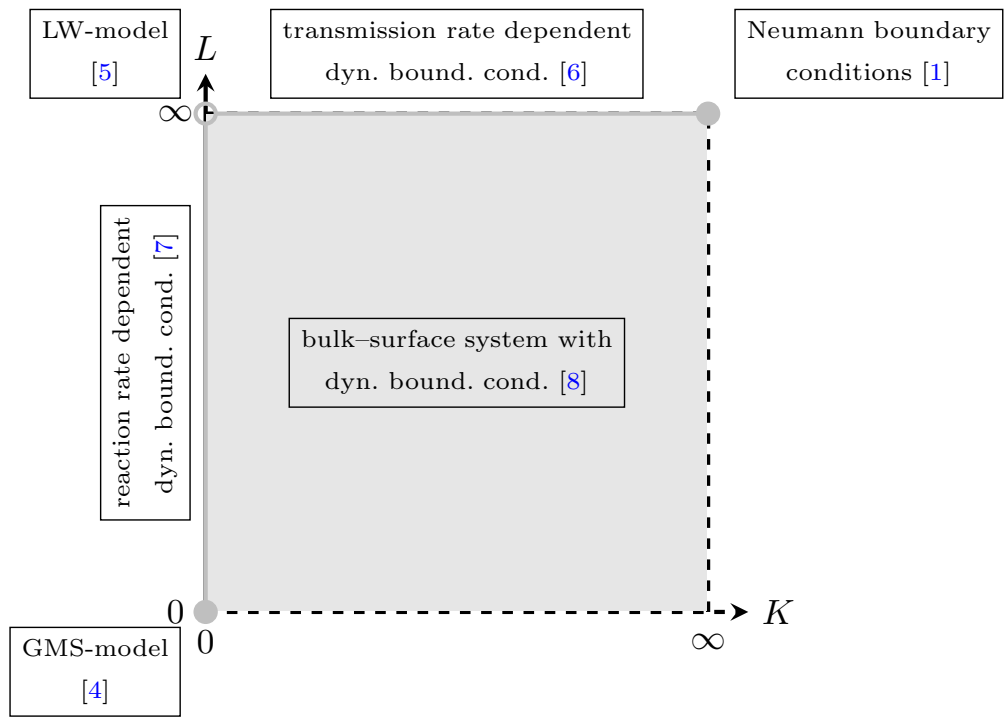

	The bulk--surface Cahn--Hilliard system with dynamic boundary conditions was proposed in \cite{KnopfStange2024} and generalizes many important systems of dynamic boundary conditions, previously developed in \cite{GoldsteinMiranvilleSchimperna2011,LiuWu2019,KnopfLam2020,KnopfLamLiuMetzger2021}, cf. Figure~\ref{fig:graphDynBC}.
	
	Let $\Om \subset \R^d$, with $d\in\{2,3\}$, and $\Ga =\partial \Om$, then we introduce the bulk phase field $u \colon \Om \times (0,T) \to \R$, the surface phase field $\psi \colon \Ga \times (0,T) \to \R$, the bulk chemical potential $\mu \colon \Om \times (0,T) \to \R$ and the surface chemical potential $\theta \colon \Ga \times (0,T) \to \R$, where the coupling of the phase fields, respectively the chemical potentials, is realized by the Robin type condition \eqref{eq:robintypecond1} and the bulk--surface coupling parameters $\alpha\in \R,K \in [0,\infty]$, and respectively by the condition \eqref{eq:robintypecond2} and the parameters $\beta\in\R,L \in [0,\infty]$. The dynamics on the boundary are of Cahn--Hilliard-type and so we arrive at the system of equations:
	\begin{subequations}\label{eq:CHreact}
		\begin{alignat}{2}
			&\partial_t u =  m_\Om \Delta \mu \qquad&&\text{ in } \Om \times (0,T), \\
			&\mu =  - \epsilon \Delta u + \epsilon^{-1} F_\Om'(u) &&\text{ in } \Om \times (0,T), \\
			&\partial_t \psi =  m_\Ga \Delta_\Ga \theta - \beta m_\Om \partial_\nu \mu  &&\text{ on } \Ga \times (0,T), \\
			&\theta =  - \delta \kappa \Delta_\Ga \psi + \delta^{-1} F_\Ga'(\psi) + \alpha \epsilon \partial_\nu u \qquad &&\text{ on } \Ga \times (0,T), \\
			&\begin{cases}
				\epsilon K\partial_\nu u =  \alpha \psi - u &\qquad \text{ if } K\in[0,\infty)\\
				\partial_\nu u = 0&\qquad \text{ if } K=\infty
			\end{cases}\qquad &&\text{ on } \Ga \times (0,T), \label{eq:robintypecond1}\\
			&\begin{cases}
				m_\Om L\partial_\nu \mu =  \beta \theta - \mu &\qquad \text{ if } L\in[0,\infty)\\
				\partial_\nu \mu = 0&\qquad \text{ if } L=\infty
			\end{cases}\qquad &&\text{ on } \Ga \times (0,T), \label{eq:robintypecond2}\\
			&(u,\psi)|_{t=0} =  (u_0,\psi_0),
		\end{alignat}
	\end{subequations}
	for initial values $u_0 \colon \Om \to \R$, $\psi_0 \colon \Ga \to \R$ and where $m_\Om, m_\Ga >0$ denote diffusion mobility parameters and $F_\Om,F_\Ga \colon \R \to \R$ are free energy potentials (e.g. double well potentials with minima at $\pm 1$). We denote by $\Delta_\Ga$ and $\partial_\nu$ the Laplace--Beltrami operator and the normal derivative on the surface respectively, as defined in \cite[Section~2.2]{DziukElliott2013}.
	
	After a short period of time, the phase field $u$ will attain values close to $\pm 1$, corresponding to the pure phases of the materials, in large regions of the domain $\Om$. These regions are separated by an interfacial region whose thickness is determined by the diffusion interface parameters $\epsilon, \delta, \kappa >0$, cf., e.g. \cite{GoldsteinMiranvilleSchimperna2011}.
	
	The Robin-type coupling \eqref{eq:robintypecond2} states that the mass flux, i.e. the motion of the materials towards and away from the boundary, is directly proportional to the difference of the bulk and surface chemical potentials. This coupling was first introduced in \cite{KnopfLamLiuMetzger2021} for the case $K=0$. The coupling of the chemical potentials is controlled by a relaxation parameter $L\in [0,\infty]$, such that $1/L$ can be interpreted as the reaction rate. By taking limits for $L$, classical dynamic boundary conditions are recovered.
	
	When considering the limit case of instantaneous reactions, so $1/L \to \infty$, one arrives at the condition that the chemical potentials are always in equilibrium:
	\begin{align*}
		\beta \theta = \mu \qquad \text{on } \Ga \times (0,T),
	\end{align*}
	hence for $\beta=1$ and $K=0$ the case of the GMS model\footnote{The abbreviation represents the surnames of the respective authors.\label{footnote:abbreviation models}}, i.e. the dynamic boundary conditions proposed in \cite{GoldsteinMiranvilleSchimperna2011}. On the other hand a vanishing reaction rate, so $1/L \to 0$, corresponds to a vanishing mass flux
	\begin{align*}
		\partial_\nu \mu = 0 \qquad \text{on } \Ga \times (0,T),
	\end{align*}
	which, for $K=0$, is the condition of the LW model\footref{footnote:abbreviation models}, i.e. the dynamic boundary conditions proposed in \cite{LiuWu2019}. It was established in \cite[Corollary 5.9]{KnopfLamLiuMetzger2021}, that in fact, if we suppose that $(u^L,\psi^L,\mu^L,\theta^L)$ is a solution of \eqref{eq:CHreact} for $L>0$ and $K=0$, then the formal limit $(u^\ast,\psi^\ast,\mu^\ast,\theta^\ast)$ as $L \to \infty$ is a solution to the LW model, and the formal limit $(u_\ast,\psi_\ast,\mu_\ast,\theta_\ast)$ as $L \to 0$ is a solution to the GMS model.
	
	The coupling of the bulk and surface phase fields is realized in the same way by \eqref{eq:robintypecond1}, which is controlled by a relaxation parameter $K\in [0,\infty]$, such that $1/K$ can be interpreted as a transmission rate. This corresponds to a case where, for instance, the bulk materials are transformed at the boundary by a chemical reaction. For $L=\infty$ this was first analyzed in \cite{KnopfLam2020}. 
	
	Similarly to the coupling of the chemical potential, when considering the limit case of instantaneous transmission, so $1/K \to \infty$, one arrives back at the equilibrium for the phase fields:
	\begin{align*}
		\alpha \psi = u \qquad \text{on } \Ga \times (0,T),
	\end{align*}
	hence for $\alpha=1$ and $L=\infty$ the case of the LW model. On the other hand a vanishing reaction rate, so $1/K \to 0$, corresponds to
	\begin{align*}
		\partial_\nu u = 0 \qquad \text{on } \Ga \times (0,T),
	\end{align*}
	which, for $L=\infty$, leads to a completely decoupled system and the classical homogeneous Neumann boundary conditions for the bulk phase field and chemical potential.
	
	It was established in \cite[Theorem 2.3]{KnopfLam2020}, that in fact, if we suppose that $(u^K,\psi^k,\mu^K,\theta^K)$ is a solution of \eqref{eq:CHreact} for $K>0$ and $L=\infty$, then the formal limit $(u^\ast,\psi^\ast,\mu^\ast,\theta^\ast)$ as $K \to 0$ is a solution to the LW model. 
	
	These considerations and asymptotic limits lead to the general formulation of the bulk--surface Cahn--Hilliard system with dynamic boundary conditions \eqref{eq:CHreact}, which was developed in \cite{KnopfStange2024}, even with an additional convection term. The existence of weak solutions for $K,L \in [0,\infty]$ and asymptotic limits for $K,L \in \{0,\infty\}$ are shown therein for a nonlinear potential with polynomial growth, see \hyperref[eq:Preliminary3]{(P3)}. The overall picture of these dynamic boundary conditions is illustrated in Figure~\ref{fig:graphDynBC}. 
	
	\section{Notation and Preliminaries}
	
	\label{section:Notation}
	\subsection*{Notation} We use the following notation throughout the paper: Let $\Om \subset \R^d$, for $d \in \{2,3\}$, with $C^2$-boundary $\Ga =\partial \Om$, then we denote the standard Lebesgue spaces for $1 \leq p \leq \infty$ and $k\geq 0$ by $L^p(\Om)$, and the Sobolev spaces by $W^{k,p}(\Om)$ (see, e.g., \cite[Section~5.2.2]{Evans1998}). For $p=2$ we use the standard notation for the Hilbert spaces $H^k(\Om)=W^{k,2}(\Om)$. For the Lebesgue and Sobolev spaces on $\Ga$ we use the similar notation $L^p(\Ga)$ and $W^{k,p}(\Ga)$ (see, e.g., \cite[Definition 2.11]{DziukElliott2013}). For any Banach space $X$, we denote by $\| \cdot \|_X$ its norm, by $X'$ the dual space and the associated duality pairing by $\langle \cdot , \cdot \rangle_X$. For a Hilbert space $X$ we denote the inner product by $(\cdot , \cdot )_X$ and if $u \in L^1(\Om)$, respectively $u \in L^1(\Ga)$, we define the spatial mean on $\Om$, respectively $\Ga$, as
	\begin{align*}
		\langle u \rangle_\Om := \frac{1}{|\Om|} \int_\Om u , \quad \text{respectively } \langle u \rangle_\Ga := \frac{1}{|\Ga|} \int_\Ga u .
	\end{align*}
	
	\subsection*{Preliminaries}
	\begin{itemize}
		\item[(P1)] \label{eq:Preliminary1}Let $\nu$ denote the unit outward pointing normal vector to $\Ga$, then we define the surface gradient $\nabla_\Ga$ of a function $u:\Ga \to \mathbb{R}$ by $\nabla_\Ga u := \nabla \bar u - (\nabla \bar u \cdot \nu)\nu$, where $\bar u$ is an arbitrary extension of $u$ in a neighborhood of $\Ga$. The Laplace--Beltrami operator on $\Ga$ is given by $\Delta_\Ga u = \nabla_\Ga \cdot \nabla_\Ga u$. Moreover, $\gamma u$ denotes the trace of $u$ on $\Ga$, for brevity we will mostly suppress the trace operator in the notation below. The normal derivative of $u$ on $\Ga$ is denoted by $\partial_\nu u$. For a more explicit definition of these operators, see \cite[Section~2.2]{DziukElliott2013} and \cite[Section~5.5]{Evans1998}.

		\item[(P2)] \label{eq:Preliminary2} We introduce for a parameter $J\in [0,\infty]$ and $\lambda \in \R$ the combined Hilbert spaces
		\begin{align*}
			&H:= L^2(\Om) \times L^2(\Ga), \quad D^\lambda:=\{ (\phi,\psi) \in H^1(\Om) \times H^1(\Ga) \colon \phi=\lambda \psi \text{ a.e. on } \Ga \}, \\
			&\text{and }\qquad  V^{J,\lambda}:=\begin{cases}
				H^1(\Om) \times H^1(\Ga), \quad &\text{if } J\in (0,\infty] \\
				D^\lambda, \quad &\text{if } J=0
			\end{cases}.
		\end{align*}
		Following the notation from \cite[Section~2.1.1]{KovacsLubich2017}, we introduce the following bilinear forms 
		\begin{equation*}
			\begin{aligned}
				m:=(\cdot,\cdot)_H \colon &H \times H \to \mathbb{R} ,\\
				& m\big((\phi,\psi),(\zeta,\xi)\big) :=  \int_\Om \phi  \zeta + \int_\Ga  \psi \xi  , \\
				a^\ast:=(\cdot,\cdot)_{V^{J,\lambda}} \colon &V^{J,\lambda} \times V^{J,\lambda} \to \mathbb{R} ,\\
				& a^\ast\big((\phi,\psi),(\zeta,\xi)\big) :=  \int_\Om \nabla \phi \cdot \nabla \zeta + \int_\Om \phi  \zeta + \int_\Ga \nabla_\Ga \psi \cdot \nabla_\Ga \xi + \int_\Ga  \psi \xi  , \\
				a^{J,\lambda} \colon &V^{J,\lambda} \times V^{J,\lambda} \to \mathbb{R} ,\\
				&a^{J,\lambda}\big((\phi,\psi),(\zeta,\xi)\big) :=  \int_\Om \nabla \phi \cdot \nabla \zeta + \int_\Ga \nabla_\Ga \psi \cdot \nabla_\Ga \xi  + \h(J) \int_\Ga (\lambda \psi-\phi)(\lambda \xi-\zeta) ,  
			\end{aligned}
		\end{equation*}
		and the semi-norm
		\begin{align*}
			|(\phi,\psi)|_{a^{J,\lambda}} :=  \Big(a^{J,\lambda}((\phi,\psi),(\phi,\psi))\Big)^{\frac12},
		\end{align*}
		where we set 
		\begin{align*}
			\h(J):= \begin{cases}
				J^{-1}, \qquad &\text{if }J\in(0,\infty),\\
				0, \qquad &\text{if }J\in \{0,\infty\}
			\end{cases}.
		\end{align*}
		
		We additionally define the weak dual norms, for $d \in H$:
		\begin{equation} \label{eq:dualNormDef} 
			\begin{aligned}
				\| d \|_{\ast,J,\lambda}:= \sup_{0\not=\varphi \in V^{J,\lambda}}\frac{m(d,\varphi)}{\|\varphi\|_{V^{J,\lambda}}}.
			\end{aligned}
		\end{equation}
		\item[(P3)] \label{eq:Preliminary3} Since the choice of $m_\Om$, $m_\Ga$, $\epsilon, \delta$ and $\kappa$ has no impact on the mathematical analysis, we simply set, without loss of generality, $m_\Om=m_\Ga=\epsilon=\delta=\kappa=1$ for the remainder of the theoretical part of this paper. Note that this does not impact any time-step restrictions, as the time discretization is linearly implicit.
		
		Following \cite{KnopfStange2024}, a weak solution of the system \eqref{eq:CHreact} on $[0,T]$, for $K,L \in [0,\infty]$ and $(u_0,\psi_0) \in V^{K,\alpha}$, is a quadruplet $(u,\psi,\mu,\theta)$, satisfying the following properties:
		\begin{enumerate}
			\item[i)] The functions $u$, $\psi$, $\mu$ and $\theta$ have the regularity:
			\begin{equation} \label{eq:WeakFormReg}
				\begin{aligned}
					&(u,\psi) \in C([0,T];H) \cap H^1(0,T;(V^{L,\beta})') \cap L^\infty(0,T;V^{K,\alpha}), \\
					&(\mu,\theta) \in L^2(0,T;V^{L,\beta}). 
				\end{aligned} 
			\end{equation} 
			\item[ii)] The functions $u$ and $\psi$ satisfy the initial conditions $u|_{t=0}=u_0$ a.e. in $\Om$, and $\psi|_{t=0}=\psi_0$ a.e. on $\Ga$.
			\item[iii)] The functions $u$, $\psi$, $\mu$ and $\theta$ satisfy the weak formulation
			\begin{subequations} \label{eq:CHweak}
				\begin{align}
					m\big((\partial_t u, \partial_t \psi),(\zeta,\xi)\big) + a^{L,\beta}\big((\mu,\theta),(\zeta,\xi)\big) =&\ 0, \label{eq:CHweak1}\\
					m\big((\mu, \theta),(\eta,\chi)\big) - a^{K,\alpha}\big((u,\psi),(\eta,\chi)\big) =&\  m\big((F_\Om'(u),F_\Ga'(\psi)),(\eta,\chi)\big),\label{eq:CHweak2}
				\end{align}
			\end{subequations}
			a.e. in $[0,T]$, for all $(\zeta,\xi) \in V^{L,\beta}$, $(\eta,\chi) \in V^{K,\alpha}$.
			\item[iv)] The functions $u$ and $\psi$ satisfy the mass conservation law
			\begin{align} \label{eq:MassConservation}
				\begin{cases}
					\beta \int_\Om u(t) + \int_\Ga \psi(t) = \beta \int_\Om u_0 + \int_\Ga \psi_0, \qquad \text{if } L\in [0,\infty), \\
					\int_\Om u(t) = \int_\Om u_0 \text{ and } \int_\Ga \psi(t) = \int_\Ga \psi_0, \qquad \text{if } L=\infty,
				\end{cases}
			\end{align}
			for all $t\in [0,T]$. This can be obtained by testing \eqref{eq:CHweak1} with $(\beta,1) \in V^{L,\beta}$ and integrating over $t$.
			\item[v)] The functions $u$ and $\psi$ satisfy the energy inequality
			\begin{align*}
				E_K(u(t),\psi(t)) + \int_0^t \int_\Om |\nabla \mu|^2 + \int_0^t \int_\Ga |\nabla_\Ga \theta|^2 + \h(L) \int_0^t \int_\Ga |\beta\theta - \mu|^2 \leq E_K(u_0,\psi_0),
			\end{align*}
			for all $t \in [0,T]$ and the Ginzburg--Landau bulk--surface free energy
			\begin{align*}
				E_K(u,\psi):=&\ \int_\Om \frac\epsilon2 |\nabla u|^2 + \int_\Om \epsilon^{-1}F_\Om(u) + \int_\Ga \frac{\delta\kappa}{2} |\nabla_\Ga \psi|^2 + \int_\Ga \delta^{-1}F_\Ga(\psi) \nonumber \\
				&\ + \h(K) \int_\Ga \frac12 |\alpha \psi - u|^2.
			\end{align*}
			This can be obtained by testing \eqref{eq:CHweak1} with $(\mu,\theta) \in V^{L,\beta}$, subtracting \eqref{eq:CHweak2} tested with $(\partial_t u,\partial_t \psi) \in V^{K,\alpha}$ and integrating over $t$.
		\end{enumerate}
		
		The main result in \cite[Theorem 3.2]{KnopfStange2024} provides existence and high regularity of weak solutions, satisfying i)--v), for the bulk--surface Cahn--Hilliard system with dynamic boundary conditions \eqref{eq:CHreact}. This result requires polynomial growth conditions on the nonlinear potentials, which include, for example, the important double well potential, and provides, for the domain $\Om$ of class $C^k$, the regularity
		\begin{equation} \label{eq:Higher Regularity}
			\begin{aligned}
				&(\mu,\theta) \in L^4(0,T;H), \quad &&\text{if } K \in (0,\infty], \\
				&(u,\psi) \in L^4(0,T;H^2(\Om) \times H^2(\Ga)), \quad &&\text{if } K \in (0,\infty], \\
				&(u,\psi) \in L^2(0,T;H^k(\Om) \times H^k(\Ga)), \quad &&\text{for } k=2,3, \\
				&(u,\psi) \in C([0,T];H^1(\Om) \times H^1(\Ga)), \quad &&\text{if } L \in (0,\infty] \text{ and } k=3.
			\end{aligned}
		\end{equation}
		\begin{rem}
			We mostly work with this abstract formulation of the problem to showcase that the following arguments can be applied in a more general setting, we therefore expect that this strategy of proof can be adapted to a wider class of problems, see, e.g., Section~\ref{section:Generalization} for an outline of the optimal-order convergence result for the Cahn--Hilliard equation on an evolving surface.
		\end{rem}
	\end{itemize}
	
	\subsection*{Assumptions}
	Let $\Om \subset \mathbb{R}^d$, for $d=2,3$, be a bounded domain with an (at least) $C^2$ boundary $\Ga=\partial \Om$.
	\begin{itemize}
		\item[(A1)] \label{eq:AssumptionPara} We assume for the parameters
		\begin{itemize}
			\item $K\in [0,\infty)$, $L\in (0,\infty)$: $\alpha\beta|\Om| + |\Ga| \not= 0$, required for the bulk--surface Poincar\'e--Wirtinger inequality, cf. Lemma~\ref{Lemma:BulkSurfPoincare}.
			\item $K\in (0,\infty)$, $L=\infty$: $\alpha \not=0$ and a $C^3$-boundary $\Ga$.
			\item $K =  0$, $L=0$: $\alpha\beta|\Om| + |\Ga| \not= 0$ and $F_\Om'(\alpha \cdot)=\beta F_\Ga'(\cdot)$ for the nonlinear potentials.
		\end{itemize}
		\item[(A2)]  \label{eq:AssumptionPot} The nonlinear potentials and their derivatives
		\begin{equation*}
			F_\Om^{(k)} \text{ and } F_\Ga^{(k)} \text{ are locally Lipschitz continuous for }k=0,\dots,3. 
		\end{equation*}
		An important example which fits into this framework is the quartic double well potential ($F(s)=\frac14 (s^2-1)^2$).
	\end{itemize}

	\section{Full discretization of the bulk--surface Cahn--Hilliard system with dynamic boundary conditions}
	\label{section:Discretization CH}
	
	\subsection{Bulk--surface finite elements} \label{sec:bulk--surfaceFEM}
	In order to discretize the equation in space, we use the linear finite element method for partial differential equations developed by Elliott and Ranner \cite[Sections 4,5]{ElliottRanner2013} and by \cite[Section~4]{KovacsLubich2017}, see also \cite[Sections 4--7]{ElliottRanner2021}. 
	
	We approximate the domain $\Om$ by a polyhedral domain $\Om_h$ and take a quasi-uniform triangulation $\calT_h$ of $\Om_h$ (see \cite[Definition 4.4.13]{BrennerScott2008}). This triangulation consists of closed simplices, either triangles in $\mathbb{R}^2$ or tetrahedra in $\mathbb{R}^3$. Let $\Ga_h:=\partial \Om_h$, then the faces of $\Ga_h$ are $(d-1)$-simplices, hence edges in $\mathbb{R}^2$ or triangles in $\mathbb{R}^3$, and we assume that the vertices lie on $\Ga$. For each point $x$ in a small, tubular neighborhood of $\Ga$ we can define the \emph{closest point projection} $p(x) \in \Ga$, cf. \cite{DziukElliott2013}. We assume the mesh width $h=\max\{\text{diam}(T) \colon T \in \calT_h\}$ to be sufficiently small, such that $\Ga_h$ is contained in this small neighborhood, then for any $x \in \Ga_h$, there exists the unique \emph{closest point} $p(x) \in \Ga$, such that $x-p(x)$ is orthogonal to the tangent space of $\Ga$ at $p$. For more details on the construction see \cite[Section~4]{ElliottRanner2013}.
	
	We define the non-conforming finite element space $V_h\nsubseteq H^1(\Om)$, corresponding to $\calT_h$, by the continuous, piecewise linear nodal basis functions on $\Om_h$: Let $(x_k)_{k=1}^N$ be the collection of all nodes of the triangulation in $\Om_h$, then we define $(\phi_j)_{j=1}^N$ as the family of continuous functions, which are piecewise linear on each of the simplices of $\calT_h$, subject to the condition
	\begin{align*}
		\phi_j(x_k)=\delta_{jk} \text{ for } j,k=1,...,N. 
	\end{align*}
	Then the \emph{bulk} finite element space is defined as $V_h:= \Span\{\phi_1,...,\phi_N\}$.
	
	If we use, without loss of generality, that the nodes of the triangulation $\calT_h$ on the boundary $\Ga_h$ are the nodes $(x_k)_{k=1}^M$, then the \emph{surface} finite element space $S_h\nsubseteq H^1(\Ga)$ is given by
	\begin{align*}
		S_h:=  \Span\{\phi_1|_{\Ga_h},...,\phi_M|_{\Ga_h}\}.
	\end{align*}
	These definitions are carried out in more detail in \cite[Section~6.2]{ElliottRanner2021}.
	
	Next we define the lift operator $\cdot^\ell \colon V_h \to V_h^\ell$ of functions $v_h \colon \Om_h \to \mathbb{R}$ to $v_h^\ell \colon \Om \to \mathbb{R}$ by setting $v_h^\ell(p)=v_h(x)$ for $x\in \Om_h$ and $p \in \Om$, which are related as described in detail in \cite[Section~4]{ElliottRanner2013}, and by defining
	\begin{align*}
		V_h^\ell:=\{v_h^\ell : v_h \in V_h\}.
	\end{align*}
	For functions $v_h \in S_h$ on the discrete surface this lift operator is then given by
	\begin{align*}
		v_h^\ell \colon \Ga \to \mathbb{R} \qquad \text{with } v_h^\ell(p)=v_h(x) ,\quad \forall p\in \Ga,
	\end{align*}
	where $x \in \Ga_h$ is the unique point on $\Ga_h$ with $x-p$ orthogonal to the tangent space $T_p\Ga$.
	
	The discrete tangential gradient $\nabla_{\Ga_h}$ is piecewisely defined analogously to $\nabla_\Ga$, see \cite[Section~4]{DziukElliott2013}. The discrete trace operator $\gamma_h v_h$ is defined by the restriction of the continuous functions $v_h \in V_h$ onto $\Ga_h$, with the same notational conventions as before, e.g. $\nabla_{\Ga_h}v_h:=\nabla_{\Ga_h}(\gamma_h v_h)$.
	
	For the combined system we define the \emph{bulk--surface} finite element space, analogous to those of the non-discretized problem in \hyperref[eq:Preliminary2]{(P2)}, for $J\in [0,\infty]$ and $\lambda \in \R$ as 
	\begin{align}
		V_{h}^{J,\lambda}:=\begin{cases}
			V_h \times S_h, \qquad &\text{if } J\in (0,\infty],\\
			D_{h}^{\lambda}:=\{(v_h,s_h) \in V_h \times S_h \colon v_h|_{\Ga_h}=\lambda s_h\}, \qquad &\text{if } J=0.
		\end{cases}\label{eq:DefFESpace}
	\end{align}

	\subsection{Linearly implicit backward difference time discretization\protect\footnote{Up to some changes, the text of this preparatory subsection is taken from \cite[Section~4]{BullerjahnKovacs2024}.}} \label{sec:BDF}
	We combine the discretization in space by the bulk--surface finite elements with the following discretizations in time: Let $q=1,\dotsc,5$, $\tau >0$ be the time step size, and $q\tau \leq t_n:=n\tau \leq T$ be a uniform partition of the time interval $[0,T]$. We assume the starting values $(u_h^0,\psi_h^0), \dots, (u_h^{q-1},\psi_h^{q-1}) \in V_h^{K,\alpha}$ to be given, and define, for $q\leq n$, the discretized time derivative, and the extrapolation for the non-linear term, respectively, as
	\begin{equation*}\label{eq:DefdiscDerivandExtra}
		\begin{aligned}
			&\partial^\tau_q u_h^n:= \frac{1}{\tau} \sum_{j=0}^q \delta_j u_h^{n-j},\qquad	&&\widetilde u_h^n:= \sum_{j=0}^{q-1} \gamma_j u_h^{n-1-j}, \\
			&\partial^\tau_q \psi_h^n:= \frac{1}{\tau} \sum_{j=0}^q \delta_j \psi_h^{n-j},\qquad	&&\widetilde \psi_h^n:= \sum_{j=0}^{q-1} \gamma_j \psi_h^{n-1-j}.
		\end{aligned}
	\end{equation*}
	Where the coefficients are given by the expressions
	\begin{subequations}\label{eq:BdfCoeff}
		\begin{align}
			\label{eq:BDF generating function}
			\delta(\zeta)=&\sum_{j=0}^q \delta_j \zeta^j=\sum_{l=1}^q \frac{1}{l} (1-\zeta)^l, \\
			\gamma(\zeta)=&\sum_{j=0}^{q-1}\gamma_j \zeta^j=\frac{(1-(1-\zeta)^q)}{\zeta}.
		\end{align}
	\end{subequations}
	We denote by $\partial^\tau:=\partial^\tau_1$ the BDF-1 derivative or backward Euler derivative, which appears several times in the following analysis independent of the choice of $q$.
	
	The $q$-step BDF methods are of order $q$, they are $A$-stable for $q=1,2$, $A(\alpha_q)$-stable for $q=3,\dotsc,6$ with $\alpha_3 =86.03^\circ$, $\alpha_4=73.35^\circ$, $\alpha_5=51.84^\circ$,  and $\alpha_6 = 17.84^\circ$, respectively, and unstable for $q \geq 7$, see, e.g.,  \cite[Section~V.2]{HairerWanner1996}, \cite[Section~2.3]{AkrivisLubich2015}, while see \cite{AkrivisKatsoprinakis2020} for the exact $\alpha_q$ values. 
	
	The starting values can be precomputed using either a lower order method with smaller time step size, or a Runge--Kutta method of suitable order.
	
	\subsection{The fully discrete scheme in abstract form}
	
	We apply the bulk--surface finite element approximation (see subsection \ref{sec:bulk--surfaceFEM}) to the abstract variational formulation \eqref{eq:CHweak} of the bulk--surface Cahn--Hilliard system with dynamic boundary conditions, namely we discretize the integrals in the bilinear forms, defined in \hyperref[eq:Preliminary2]{(P2)} for $J\in[0,\infty]$ and $\lambda \in \R$, by integrals over the polyhedral domains $\Om_h$ and $\Ga_h$, following \cite[Section~4.1]{KovacsLubich2017}, for the space $V_h^{J,\lambda}$ defined in \eqref{eq:DefFESpace}:
	\begin{equation*} \label{eq:DefDiscBilinearForms}
		\begin{aligned}
			m_h \colon &V_h^{J,\lambda} \times V_h^{J,\lambda} \to \mathbb{R} ,\\
			&m_h\big((\phi_h,\psi_h),(\zeta_h,\xi_h)\big) :=  \int_{\Om_h} \phi_h \zeta_h + \int_{\Ga_h} \psi_h \xi_h ,\\
			a_h^{\ast} \colon &(V_h\times S_h) \times (V_h \times S_h) \to \mathbb{R} ,\\
			&a_h^{\ast}\big((\phi_h,\psi_h),(\zeta_h,\xi_h)\big) := \int_{\Om_h} \nabla \phi_h \cdot \nabla \zeta_h + \int_{\Om_h} \phi_h \zeta_h +  \int_{\Ga_h} \nabla_{\Ga_h} \psi_h \cdot \nabla_{\Ga_h} \xi_h + \int_{\Ga_h} \psi_h \xi_h  , \\	
			a_h^{J,\lambda} \colon &V_h^{J,\lambda} \times V_h^{J,\lambda} \to \mathbb{R} ,\\
			&a_h^{J,\lambda}\big((\phi_h,\psi_h),(\zeta_h,\xi_h)\big) := \int_{\Om_h} \nabla \phi_h \cdot \nabla \zeta_h + \int_{\Ga_h} \nabla_{\Ga_h} \psi_h \cdot \nabla_{\Ga_h} \xi_h  \nonumber \\
			&\hphantom{\big((\phi_h,\psi_h),(\zeta_h,\xi_h)\big) \mapsto} + \h(J) \int_{\Ga_h} (\lambda \psi_h - \phi_h) (\lambda \xi_h-\zeta_h) .	 
		\end{aligned}
	\end{equation*}
	Note that we define the discrete norms and the semi-norm analogous to the continuous case by
	\begin{align}
		\| \cdot \|_{H_h}^2:=m_h(\cdot,\cdot), \quad \| \cdot \|_{V_h^{J,\lambda}}^2:=a_h^\ast(\cdot,\cdot), \quad |\cdot |_{a_h^{J,\lambda}}^2:=a_h^{J,\lambda}(\cdot,\cdot). \label{eq:DefDiscNorms}
	\end{align}
	
	Then we use a Galerkin ansatz in order to solve this variational problem in a finite element space and combine this with the discretization in time by the $q$-step linearly implicit backwards difference method (see Section \ref{sec:BDF}). Then the fully discrete scheme reads: Let $L,K \in [0,\infty]$, then find, for $h,\tau >0$ and $q \tau \leq t_n:=n\tau \leq T$, functions $(u_h^n, \psi_h^n) \in V_h^{K,\alpha}$ and $(\mu_h^n,\theta_h^n) \in V_h^{L,\beta}$ such that
	\begin{subequations} \label{eq:CHweakD}
		\begin{align}
			&m_h\big((\partial_q^\tau u_h^n, \partial_q^\tau \psi_h^n),(\zeta_h,\xi_h)\big) + a_h^{L,\beta}\big((\mu_h^n,\theta_h^n),(\zeta_h,\xi_h)\big) = 0, \label{eq:CHweakD1}\\
			&m_h\big((\mu_h^n, \theta_h^n),(\eta_h,\chi_h)\big) - a_h^{K,\alpha}\big((u_h^n,\psi_h^n),(\eta_h,\chi_h)\big) = m_h\big((\calI_{\Om_h} F_\Om'(\widetilde u_h^n),\calI_{\Ga_h} F_\Ga'(\widetilde\psi_h^n)), (\eta_h,\chi_h)\big).\label{eq:CHweakD2}
		\end{align}
	\end{subequations}
	for all $(\zeta_h,\xi_h) \in V_h^{L,\beta}$ and $(\eta_h,\chi_h) \in V_h^{K,\alpha}$ and with initial values $(u_h^0,\psi_h^0),...,(u_h^{q-1},\psi_h^{q-1}) \in V_h^{K,\alpha}$ given, and where $\calI_{\Om_h} \colon C(\Om_h) \to V_h$ and $\calI_{\Ga_h}\colon C(\Ga_h) \to S_h$ denote standard finite element interpolation operators, see e.g. \cite{BrennerScott2008}.

	\subsection{Geometric error estimates} \label{Subsection:GeomErrorEst}
	We recall, that the bounded domain $\Om \subset \mathbb{R}^d$ $(d=2,3)$ has an at least $C^2$-boundary $\Ga$, and $\Ga_h=\partial \Om_h$, the boundary of the quasi-uniform triangulation $\Om_h$, is an interpolation of $\Ga$. 
	The following estimates were shown in \cite[Lemma~6.2]{ElliottRanner2013}.
	\begin{lem}[geometric approximation errors] \label{Lemma:GeomAppError}
		The following estimates hold for $h \leq h_0$ and for any $v_h,w_h \in V_h$ and $\phi_h,\eta_h \in S_h$,
		\begin{align*}
			\Big|\int_\Om \nabla v^\ell_h \cdot \nabla w^\ell_h - \int_{\Om_h}\nabla v_h\cdot \nabla w_h\Big| \leq&\  C h \|\nabla v_h^\ell\|_{L^2(B_h^\ell)} \|\nabla w_h^\ell\|_{L^2(B_h^\ell)}, \\
			\Big|\int_\Ga \nabla_\Ga \phi^\ell_h\cdot \nabla_\Ga \eta^\ell_h - \int_{\Ga_h}\nabla_{\Ga_h}\phi_h\cdot \nabla_{\Ga_h} \eta_h \Big| \leq&\ C h^2 \| \nabla_\Ga \phi_h^\ell \|_{L^2(\Ga)} \| \nabla_\Ga \eta_h^\ell \|_{L^2(\Ga)}, \\
			\Big|\int_\Om v^\ell_h w^\ell_h  - \int_{\Om_h} v_h w_h \Big| \leq&\   C h \|v_h^\ell\|_{L^2(B^\ell_h)} \| w_h^\ell\|_{L^2(B^\ell_h)}, \\
			\Big|\int_\Ga \phi^\ell_h \eta^\ell_h  - \int_{\Ga_h} \phi_h \eta_h \Big| \leq&\  Ch^2 \|  \phi^\ell_h \|_{L^2(\Ga)} \|  \eta^\ell_h \|_{L^2(\Ga)},
		\end{align*}
		where $B^\ell_h$ collects the lifts of elements which have at least two nodes on the boundary. 
		
		A combination of these estimates yields geometric approximation estimates between the bilinear forms in \hyperref[eq:Preliminary2]{(P2)} and their discrete counterparts in \eqref{eq:DefDiscBilinearForms}, with the appropriate order.
	\end{lem}
	
	\begin{rem} \label{rem:normequivalence}
		As a consequence we also have the $h$-uniform equivalence of the norms $\| (\cdot)^\ell \|_{H^1(\Om)}$ and $\| \cdot \|_{H^1(\Om_h)}$, and of the norms $\| (\cdot)^\ell \|_{L^2(\Om)}$ and $\| \cdot \|_{L^2(\Om_h)}$ on $V_h$, and also on the surface of the norms $\| (\cdot)^\ell \|_{H^1(\Ga)}$ and $\| \cdot \|_{H^1(\Ga_h)}$, and of the norms $\| (\cdot)^\ell \|_{L^2(\Ga)}$ and $\| \cdot \|_{L^2(\Ga_h)}$ on $S_h$.
	\end{rem}

	\subsection{Temporal error estimates \protect\footnote{Up to some changes, the text of this preparatory subsection is taken from \cite[Section~7.2]{BullerjahnKovacs2024}.}} \label{sec:TempErrorEst}
	Recalling the definition of the discrete derivative and the extrapolation, see Section~\ref{sec:BDF}, with the coefficients \eqref{eq:BdfCoeff}, we can define this in a more general way for a function $y:[0,T] \to \mathbb{R}$ as 
	\begin{align*}
		\partial^\tau_q y(s)= \frac{1}{\tau} \sum_{i=0}^q \delta_i y(s-i\tau), \qquad \widetilde y (s) = \sum_{j=0}^{q-1} \gamma_j y(s-j\tau-\tau),
	\end{align*}
	for $s\geq q\tau$, which is consistent with the previous definition for $s=n\tau$ and the natural choice $u_h^j=y(j\tau)$, for $j=n-q,...,n$. 
	
	The main tool in analyzing the temporal error in the consistency analysis is the following result, based on the Peano kernel of a multi-step method, which can be found in, e.g., \cite[Section~III.2]{HairerNorsettWanner1993}.
	\begin{lem}
		\label{Lemma:PeanoKernel}
		(i) Let $y\colon[0,T] \to \mathbb{R}$ be $p$ times continuously differentiable, with $y^{(p)}$ absolutely continuous, for $1 \leq p \leq q \leq 5$, and let $t^\ast \in [q\tau,T]$, then 
		\begin{align*}
			\partial^\tau_q y(t^\ast) - \frac{\d}{\d t}y(t^\ast)=\tau^p \int_{0}^{q} K_q(s) y^{(p+1)}(t^\ast-s\tau) \d s ,
		\end{align*}
		where
		\begin{align*}
			K_{q}(s)= - \frac{1}{p!} \sum_{i=0}^q \delta_i (s-i)_-^p , \qquad \text{ with } \qquad (x)_- := \begin{cases}
				x, & x \leq 0 , \\
				0, & x \geq 0 .
			\end{cases}
		\end{align*}
		(ii) Let $y\colon[0,T] \to \mathbb{R}$ be $p-1$ times continuously differentiable, with $y^{(p-1)}$ absolutely continuous, for $1 \leq p \leq q \leq 5$, and let $t^\ast \in [q\tau,T]$, then 
		\begin{align*}
			\widetilde y(t^\ast) - y(t^\ast)=\tau^p \int_{0}^{q} \widetilde K_q(s) y^{(p)}(t^\ast-s\tau) \d s,
		\end{align*}
		where 
		\begin{align*}
			\widetilde K_q(s) = \sum_{i=0}^{q-1} \gamma_i \frac{(s-(i+1))_-^{p-1}}{(p-1)!}.
		\end{align*}
	\end{lem}

	\section{Main result: Optimal-order fully discrete error estimates}
	\label{section:MainResult}
	The main result of this paper establishes optimal-order error estimates for the bulk--surface finite element and linearly implicit backward difference discretization of the bulk--surface Cahn--Hilliard system with dynamic boundary conditions \eqref{eq:CHweak}, for the range of parameters depicted in Figure~\ref{fig:graphDynBC}, cf. \hyperref[eq:AssumptionPara]{(A1)}. For the GMS-model, i.e. the choice $K=0$ and $L=0$, this result was already obtained in \cite{BullerjahnKovacs2024}, by a different technique of exploiting the anti-symmetric structure of the equation, which is not possible in the general case. To our knowledge there are no fully discrete error estimates available for other choices of parameters $K$ and $L$, note however that convergence of sub-sequences without rates are available, see, e.g. \cite{KnopfLamLiuMetzger2021}.
	
	\begin{thm}\label{thm:OptimalOrderErrorEst}
		Let $1\leq q \leq 5$ and $(u,\psi,\mu,\theta)$ be a sufficiently smooth solution (see \eqref{eq:regularityAss}) of the bulk--surface Cahn--Hilliard system with dynamic boundary conditions \eqref{eq:CHweak}, with \hyperref[eq:AssumptionPara]{(A1)} and nonlinear potentials satisfying \hyperref[eq:AssumptionPot]{(A2)}. Then there exists $h_0>0$ such that for all $h\leq h_0$ and $\tau >0$, satisfying the mild Courant--Friedrichs--Lewy (CFL) condition $\tau^q\leq C_0 h^2$ (where $C_0$ can be chosen arbitrarily), the error between the solution $(u,\psi,\mu,\theta)$ and the fully discrete solution $(u_h^n,\psi_h^n,\mu_h^n,\theta_h^n)$ of \eqref{eq:CHweakD}, using linear bulk--surface finite element discretization and linearly implicit backward difference time discretization of order $q=1,\dots,5$, satisfies the optimal-order error estimates, for $q\tau \leq t_n:=n\tau \leq T$,
		\begin{align*}
			\| &(u^n_h)^\ell - u(t_n) \|_{L^2(\Om)} + \| (\psi^n_h)^\ell - \psi(t_n) \|_{L^2(\Ga)} \nonumber \\
			&+ h \Big( \| (u^n_h)^\ell - u(t_n) \|_{H^1(\Om)} + \|  (\psi^n_h)^\ell - \psi(t_n) \|_{H^1(\Ga)}\Big) \leq C (h^2+\tau^q),  \\
			\bigg( \tau \sum_{k=q}^n& \| (\mu^k_h)^\ell - \mu(t_k) \|_{L^2(\Om)}^2 + \|(\theta^k_h)^\ell - \theta(t_k) \|_{L^2(\Ga)}^2 \nonumber \\
			&+ h \Big(\| (\mu^k_h)^\ell - \mu(t_k) \|_{H^1(\Om)}^2 + \| (\theta^k_h)^\ell - \theta(t_k) \|_{H^1(\Ga)}^2\Big) \bigg)^{\frac{1}{2}} \leq C (h^2+\tau^q),
		\end{align*} 
		provided the error in the starting values $\calE_I^{q-1}$, defined in \eqref{eq:DefIhq}, is $O((\tau^q+h^2)^2)$. The constant $C>0$ depends on the Sobolev norms of the exact solution and exponentially on the final time $T$, but is independent of $h$ and $\tau$.
	\end{thm}
	
	\begin{rem}
		(i) The $h_0$ depends on geometric approximation results (see, e.g., Section~\ref{sec:bulk--surfaceFEM} and Section~\ref{section:Consistency}). The requirement on the smallness of both parameters is merely technical. 
		
		(ii) The mild CFL condition is needed to establish an $L^\infty$-bound on the error in the stability analysis. This lets us reduce the assumptions on the nonlinear potentials $F_\Om$ and $F_\Ga$ to a \emph{local} Lipschitz condition, which covers important examples such as the double well potential. The restriction can be even weakened to $\tau^q \leq C_0 h^{3/2+\epsilon_0}$, where $\epsilon_0,C_0 > 0$ can be chosen arbitrarily, as is apparent from the assumptions in Proposition~\ref{prop:stability}.
		
		(iii) The constant $C$ also depends on a polynomial expression of the inverses of $\epsilon$ and $\delta$, that is, sharp interface limits are not covered by our error estimates.
		
		(iv) Note that the sufficient regularity assumptions \eqref{eq:regularityAss} for the exact solution are not covered by the theoretically established regularity in \cite[Theorem~3.2]{KnopfStange2024}, cf. \eqref{eq:Higher Regularity}, even though the regularity in space for the phase fields is achieved.
	\end{rem}
	
	Sufficient regularity assumptions for Theorem~\ref{thm:OptimalOrderErrorEst} are: The exact solutions satisfy the required regularity for the weak formulation \eqref{eq:WeakFormReg} and additionally
	\begin{equation} \label{eq:regularityAss}
		\begin{alignedat}{3}
			&\ u \in H^{q+2}([0,T];H^2(\Om)), & \quad &
			\psi \in H^{q+2}([0,T];H^2(\Ga)), \\
			&\ \mu \in H^{3}([0,T];H^2(\Om)), & \quad &
			\theta \in H^{3}([0,T];H^2(\Ga)).
		\end{alignedat}
	\end{equation} 
	Recall that by standard theory we have, for any $u \in H^1([0,T];X)$, that $u \in C([0,T],X)$ with
	\begin{align*} 
		\max_{0\leq t \leq T} \| u(t)\|_X \leq c(T) \|u\|_{H^1([0,T];X)},
	\end{align*}
	see, e.g., \cite[Section~5.9.2]{Evans1998}. Since the Sobolev embedding $H^2 \subset L^\infty$ is continuous, by $d=2,3$ (see, e.g., \cite[Section~5.6.3]{Evans1998}), we have the regularity
	\begin{equation*}
		\begin{aligned}
			u &\in W^{q+1,\infty}([0,T];L^{\infty}(\Om)), \quad && \psi \in  W^{q+1,\infty}([0,T];L^{\infty}(\Ga)), \\
			\text{and} \qquad \mu &\in W^{2,\infty}([0,T];L^\infty(\Om)), \quad &&~~\theta \in  W^{2,\infty}([0,T];L^\infty(\Ga)).
		\end{aligned}
	\end{equation*}

		%
		%
		%
		%
		%
	
	\refstepcounter{scheme}
	\begin{figure}[ht]	
		\renewcommand{\figurename}{Scheme} 
		The outline of our argument is as follows: 
		\begin{itemize} 
			\item[\textbf{Part A}:] Derive an energy equality by testing the error equations with the errors and the discrete derivative of the errors
			
			
			\item[\textbf{Part B}:] Use the technique based on the $G$-stability theory by Dahlquist \cite[Theorem~3.3]{Dahlquist1978} and the multiplier techniques of Nevanlinna and Odeh \cite[Section~2]{NevanlinnaOdeh1981}, sum over all time-steps and get the $H^1$-seminorm, of the error at time $t_n$ on the left-hand side of the inequality
			
			\item[\textbf{Part C}:] Eliminate terms in the discrete derivative of errors by a discrete product rule and re-substituting the equation, and then estimating the remaining terms
			
			\item[\textbf{Part D}:] Derive a Poincar\'e-type estimate by an \emph{almost mass conservation} of the error, to estimate the $H^1$-seminorm by the full $H^1$-norm 
			
			\item[\textbf{Part E}:] An absorption argument and Gr\"onwall's inequality yields the final estimate
		\end{itemize}
		\caption{General scheme of the new stability proof.}
		\label{scheme:Stability-scheme}
	\end{figure}

	In the following we will prove Theorem~\ref{thm:OptimalOrderErrorEst} by a combination of a stability and consistency analysis.
	
	Scheme~\ref{scheme:Stability-scheme} captures the main new idea for the stability proof in this paper. Comparing this approach to the stability proof for the Cahn--Hilliard equation with Cahn--Hilliard type dynamic boundary conditions in \cite[Section~6]{BullerjahnKovacs2024}, we replace a set of four energy estimates by only one energy estimate. The previously used four estimates depend heavily on the anti-symmetric structure of the equation, which is no longer present in our more general case. Instead, we use only one energy estimate and show a discrete Poincar\'e-type inequality via an \emph{almost mass conservation}, i.e. the errors caused by perturbations of the weak formulation have mass bounded by norms of the perturbations and initial data, together with the bulk--surface Poincar\'e--Wirtinger inequality, Lemma~\ref{Lemma:BulkSurfPoincare}.
	
	We therefore strongly believe that the stability Scheme~\ref{scheme:Stability-scheme} translates well to other equations satisfying an energy law and \emph{almost} mass conservation properties, allowing for similar stability results in a large class of problems, see, e.g., Section~\ref{section:Generalization} for a proof of the optimal-order convergence result for the Cahn--Hilliard equation on an evolving surface. 
	
	To our knowledge this is the first time a fully discrete Poincar\'e-type inequality was used to show stability and convergence estimates for a higher order BDF time discretization. 
	
	The transfer of the energy estimates from the continuous case to the fully discrete case is done by a method combining the $G$-stability theory by Dahlquist \cite[Theorem~3.3]{Dahlquist1978} and the multiplier techniques of Nevanlinna and Odeh \cite[Section~2]{NevanlinnaOdeh1981}, which has already proved to be effective for different partial differential equations, see, e.g., \cite{LubichMansourVenkataraman2013,AkrivisLubich2015,KovacsLiLubich2019,LLG,ContriKovacsMassing2023,HMHF,BullerjahnKovacs2024}.

	In the consistency part we are concerned with the approximation error which comes from discretizations, and enters in the stability analysis as defect terms. The optimal-order fully discrete consistency bounds are obtained analogously to \cite[Section~7]{BullerjahnKovacs2024}, by error estimates for the spatial discretization, see, e.g. \cite[Section~6]{HarderKovacs2022}, and a sufficiently regular extension of the extrapolation to obtain error estimates for the backward difference time discretization.
	
	A crucial advantage is that no shift by initial data is needed in the error analysis for the fully discrete scheme \eqref{eq:CHweakD} in comparison to the error analysis in \cite{KovacsLiLubich2021,HarderKovacs2022,BeschleKovacs2022,BullerjahnKovacs2024}.

	\section{Bulk--surface Poincar\'e--Wirtinger inequality} \label{section:BS - Poincare-wirtinger}
	
	The key technical estimate in our stability proof is the Poincar\'e--Wirtinger inequality for the combined bulk and surface spaces introduced in \hyperref[eq:Preliminary2]{(P2)}. This result is the standard reformulation of the Poincar\'e-type estimate---somewhat hidden---in the appendix of \cite[Lemma~A.1]{KnopfLiu2021}, an earlier version for $J=0$ can be found in the appendix in \cite[Lemma~A]{ColliFukao2015}.
	
	We define the generalized mean for bulk and surface functions, for $\lambda,\upsilon \in \R$, as
	\begin{align*}
		M^{\lambda,\upsilon}\colon V^{J,\lambda} \to \R, \qquad (\phi,\psi) \mapsto \frac{\upsilon |\Om| \langle \phi \rangle_\Om + |\Ga| \langle \psi \rangle_\Ga}{\lambda \upsilon |\Om| + |\Ga|} ,
	\end{align*}
	and obtain the following result.
	\begin{lem}[\bf Poincar\'e--Wirtinger inequality for bulk and surface] \label{Lemma:BulkSurfPoincare}
		Let $J\in [0,\infty)$ and $\lambda,\upsilon \in \R$, with $\lambda \upsilon |\Om| + |\Ga| \not=0$, and let $\Om \subset \R^d$ be a bounded Lipschitz domain with boundary $\Ga=\partial \Om$. Then there exists a constant $c_P >0$ depending only on $K$, $\lambda$, $\upsilon$ and $\Om$ such that 
		\begin{align*} 
			\Big\|\Big(\phi-\lambda M^{\lambda,\upsilon}((\phi,\psi)),\psi-M^{\lambda,\upsilon}((\phi,\psi))\Big)\Big\|_{V^{J,\lambda}} \leq c_P |(\phi,\psi)|_{a^{J,\lambda}}
		\end{align*}
		for all $(\phi,\psi) \in V^{J,\lambda}$.
	\end{lem}
	\begin{proof}
		The original proof in \cite[Lemma~A.1]{KnopfLiu2021} suffices here up to the mean zero transformations.
	\end{proof}
	In particular, this means that we obtain a Poincar\'e--Wirtinger inequality for the \emph{combined} bulk and surface, with the same combined mass $M^{\lambda,\upsilon}((\phi,\psi))$ in the bulk and on the surface. 
	
	Note that if one would simply use the Poincar\'e-inequality in the bulk and on the surface \emph{independently}, one would need to use different masses for the bulk and surface part. This is not feasible since this problem in general only satisfies a conservation of the \emph{combined mass}, see \eqref{eq:MassConservation}.

	\section{Stability of the numerical scheme}
	\label{section:Stability}
	
	\subsection{Ritz map and error equation} \label{sec:Ritzmap}
	In order to perform a stability analysis of the system \eqref{eq:CHweakD}, we assume that the quadruple $(u,\psi,\mu,\theta)$ is an exact solution to the equation \eqref{eq:CHweak} and has sufficiently high regularity (e.g. \eqref{eq:regularityAss}). Then we want to project this exact solution onto the discrete space. In order to do this we introduce, similar to \cite[Section~4]{KovacsLubich2017}, the bilinear forms
	\begin{equation*}
		\begin{aligned}
			&a_{\Om}^\ast \colon H^1(\Om) \times H^1(\Om) \to \mathbb{R} ,\quad &&a_{\Om}^\ast(v,w) := \int_{\Om} v w + \int_{\Om} \nabla v \cdot \nabla w ,\\
			&a_{\Ga}^\ast \colon H^1(\Ga) \times H^1(\Ga) \to \mathbb{R} ,\quad &&a_{\Ga}^\ast(v,w) := \int_{\Ga} v w + \int_{\Ga} \nabla_\Ga v \cdot \nabla_\Ga w ,\\
			&a_{\Om_h}^\ast \colon V_h \times V_h \to \mathbb{R} ,\quad &&a_{\Om_h}^\ast (\phi_h,\psi_h) := \int_{\Om_h} \phi_h \psi_h + \int_{\Om_h} \nabla\phi_h \cdot \nabla\psi_h ,\\
			&a_{\Ga_h}^\ast \colon S_h \times S_h \to \mathbb{R} ,\quad &&a_{\Ga_h}^\ast(\phi_h,\psi_h):=  \int_{\Ga_h} \phi_h \psi_h + \int_{\Ga_h} \nabla_{\Ga_h}\phi_h \cdot \nabla_{\Ga_h}\psi_h ,\\		 
		\end{aligned}
	\end{equation*}
	and define the bulk Ritz map $\widetilde R_h^\Om \colon H^1(\Om) \to V_h $, for $\zeta \in H^1(\Om)$, by the identity
	\begin{align}
		a_{\Om_h}^\ast(\widetilde R_h^\Om (\zeta) , \varphi_h)= a_\Om^\ast(\zeta,\varphi_h^\ell) \qquad \text{for every } \varphi_h \in V_h, \label{eq:DefRitzmap1}
	\end{align}
	and the surface Ritz map $\widetilde R_h^\Ga \colon H^1(\Ga) \to S_h $, for $\xi \in H^1(\Ga)$, by the identity
	\begin{align}
		a_{\Ga_h}^\ast(\widetilde R_h^\Ga (\xi) , \phi_h)= a_\Ga^\ast(\xi,\phi_h^\ell) \qquad \text{for every } \phi_h \in S_h. \label{eq:DefRitzmap2}
	\end{align}
	These maps are well defined by the Riesz representation theorem, using the ellipticity of the bilinear forms (for more details see \cite[Section~4]{KovacsLubich2017}), and we denote the (lifted) Ritz maps by $R_h^\Om(\zeta):=\big(\widetilde R_h^\Om (\zeta)\big)^\ell$ and $R_h^\Ga(\xi):=\big(\widetilde R_h^\Ga (\xi)\big)^\ell$.
	
	Then we denote by $(u^n_\ast,\psi^n_\ast)=(\widetilde R_h^\Om(u(t_n)),\widetilde R_h^\Ga(\psi(t_n)))$ and $(\mu_\ast^n, \theta_\ast^n)=(\widetilde R_h^\Om(\mu(t_n)),\widetilde R_h^\Ga(\theta(t_n)))$ the Ritz projections of the exact solution at time $t_n=n\tau$. These solve the numerical scheme only up to some defects, which we denote by $(d_1^n, d_2^n), (d_3^n, d_4^n)$ corresponding to the two equations of the scheme \eqref{eq:CHweakD}. 
	
	We arrive at the following equation system for the errors $e_u^n:=u^n - u_\ast^n$, $e_\psi^n:=\psi^n - \psi_\ast^n$, $e_\mu^n:=\mu^n - \mu_\ast^n$ and $e_\theta^n:=\theta^n-\theta_\ast^n$:
	\begin{subequations} \label{eq:erroreq}
		\begin{align}
			m_h\big((\partial_q^\tau e_u^n, \partial_q^\tau e_\psi^n),(\zeta_h,\xi_h)\big) + a_h^{L,\beta}\big((e_\mu^n,e_\theta^n),(\zeta_h,\xi_h)\big) =&\ m_h\big((d_1^n, d_2^n),(\zeta_h,\xi_h)\big), \label{eq:erroreq1}\\
			m_h\big((e_\mu^n, e_\theta^n),(\eta_h,\chi_h)\big) - a_h^{K,\alpha}\big((e_u^n,e_\psi^n),(\eta_h,\chi_h)\big) =&\  m_h\big((r_{\Om_h}^n,r_{\Ga_h}^n),(\eta_h,\chi_h)\big) \nonumber \\
			&\ + m_h\big((d_3^n, d_4^n),(\eta_h,\chi_h)\big).\label{eq:erroreq2}
		\end{align}
	\end{subequations}
	for all $(\zeta_h,\xi_h) \in V_h^{L,\beta}$ and $(\eta_h,\chi_h) \in V_h^{K,\alpha}$ and where we collect the nonlinear terms in $r_{\Om_h}^n:= \calI_{\Om_h}F_\Om'(\widetilde u_h^{n})-\calI_{\Om_h}F_\Om'(\widetilde u_\ast^{n})$ and $r_{\Ga_h}^n:= \calI_{\Ga_h}F_\Ga'(\widetilde \psi_h^{n})-\calI_{\Ga_h}F_\Ga'(\widetilde \psi_\ast^{n})$. The starting values for this scheme are given by 
	\begin{align*}
		e_u^i=u_h^i - u^i_\ast \text{ for } 0\leq i \leq q-1,\\
		e_\psi^i=\psi_h^i - \psi^i_\ast \text{ for } 0\leq i \leq q-1,
	\end{align*}
	where we assume $u_h^i$ and $\psi_h^i$ are known for $0 \leq i \leq q-1$.
	
	In order to bound the defect terms we need to define the norms in which we will estimate them. We define the weak norm for $J \in [0,\infty]$, $\lambda \in \R$ and $d_h \in V_h^{J,\lambda}$, analogous to \eqref{eq:dualNormDef}, as
	
	\begin{equation} \label{eq:DefectNorms} 
		\begin{aligned}
			\| d_h \|_{\ast,J,\lambda}:= \sup_{0\not=\varphi_h \in V_h^{J,\lambda}}\frac{m_h(d_h,\varphi_h)}{\|\varphi_h\|_{V_h^{J,\lambda}}},
		\end{aligned}
	\end{equation}
	where the dependence on the parameter $h$ is suppressed in the notation.
	
	Additionally, we set for $0 \leq i \leq q-1$ the defects as the value of the defects in the corresponding semi-discrete problem at the interpolation points, i.e. $d^i_1$, $d^i_2$, $d_3^i$ and $d^i_4$ are the vectors of nodal values of $d^1_h(i\tau)$, $d^2_h(i\tau)$, $d^3_h(i\tau)$ and $d^4_h(i\tau)$.
	
	\subsection{An $L^\infty$ bound for the Ritz map} \label{subsec:Linftyu}
	
	To control the \emph{locally Lipschitz continuous} potentials $F_\Om'$ and $F_\Ga'$ of the extrapolated value of the numerical and the exact solution, we show that all inputs are in a compact set, which allows us to use Lipschitz continuity for these terms. As a preparation we show here an $h$-uniform $L^\infty$-bound of the Ritz map of the exact solutions $(u_\ast^n,\psi_\ast^n) = (\widetilde R_h^\Om(u(t_n)),\widetilde R_h^\Ga(\psi(t_n)))$.
	
	In order to derive this bound we need the following set of estimates, concerned with the error in approximation by the Ritz map.
	\begin{lem}[Ritz map errors] \label{Lemma:Ritzmaperror}
		The Ritz maps \eqref{eq:DefRitzmap1}, and \eqref{eq:DefRitzmap2}, satisfy the following error bounds, for $h \leq h_0$:
		\begin{enumerate}
			\item \cite[Lemma~3.8]{ElliottRanner2021} For any $v \in H^2(\Om)$,
			\begin{align*}
				\| v- R_h^\Om v\|_{L^2(\Om)} + h \| v-R_h^\Om v \|_{H^1(\Om)} \leq Ch^2 \| v \|_{H^2(\Om)}.
			\end{align*}
			\item \cite[Theorem~3.2]{ElliottRanner2015} For any $v \in H^2(\Ga)$,
			\begin{align*}
				\| v- R_h^\Ga v\|_{L^2(\Ga)} + h \| v-R_h^\Ga v \|_{H^1(\Ga)} \leq Ch^2 \| \gamma v \|_{H^2(\Ga)}.
			\end{align*}
		\end{enumerate}
		where the constants $C$ are independent of $h$ and $v$.
	\end{lem}

	By the inverse estimate (see, e.g., \cite[Theorem~4.5.11]{BrennerScott2008}), the $h$-uniform norm equivalence between the discrete and original $L^2$-norms (see Remark~\ref{rem:normequivalence}), the $L^2$-norm error estimates of the finite element interpolation operator $\widetilde I_h v \in V_h$, with lift $I_h v =(\widetilde I_h v)^\ell \in V_h^\ell$ (\cite[Proposition~2.7]{Demlow2009}), and the error estimates of the Ritz map (Lemma~\ref{Lemma:Ritzmaperror}) we estimate, for any $0 \leq t \leq T$ (suppressed in the notation) and $j = 0, 1$,
	\begin{align*}
		\| \widetilde R_h^\Om u^{(j)} - \widetilde I_h u^{(j)} \|_{L^\infty(\Om_h)} \leq&\  c h^{-d/2} \| \widetilde R_h^\Om u^{(j)} - \widetilde I_h u^{(j)} \|_{L^2(\Om_h)} \\
		\leq&\  c h^{-d/2} \| R_h^\Om u^{(j)} - u^{(j)} \|_{L^2(\Om)} 
		+ c h^{-d/2} \| u^{(j)} - I_h u^{(j)} \|_{L^2(\Om)} \\
		\leq&\ c h^{2-d/2} \| u^{(j)} \|_{H^2(\Om)},
	\end{align*}
	recall that $d=2,3$ denotes the dimension of the domain $\Om \subset \R^d$, cf. \hyperref[eq:AssumptionPara]{(A1)}.
	
	Using the regularity of the exact solution \eqref{eq:regularityAss}, and the $L^\infty$-stability of the linear finite element interpolation operator, for any $0 \leq t \leq T$ and $j = 0, 1$, we have
	\begin{equation}
		\label{eq:exactuLinfty}
		\begin{aligned}
			\| u_\ast^{(j)} \|_{L^\infty(\Om_h)} \leq &\ \| \widetilde R_h^\Om u^{(j)} - \widetilde I_h u^{(j)} \|_{L^\infty(\Om_h)} 
			+ c \| I_h u^{(j)} \|_{L^\infty(\Om_h)} \\
			\leq&\  c h^{2-d/2} \| u^{(j)} \|_{H^2(\Om)} + c \| u^{(j)} \|_{L^\infty(\Om)} \\
			\leq&\ c.
		\end{aligned}
	\end{equation}
	We note that for the extrapolation, see Section~\ref{sec:BDF}, of the Ritz map of the exact solution, with $k\geq q$, we simply have
	\begin{align*}
		\| \widetilde u^k_\ast \|_{L^\infty(\Om_h)} \leq \sum_{j=0}^{q-1} |\gamma_j|  \|u_\ast^{k-1-j}\|_{L^\infty(\Om_h)} \leq c.
	\end{align*}
	
	Following these lines we establish the similar estimate,
	\begin{align}
		\| \nabla \widetilde u^k_\ast \|_{L^p(\Om_h)} \leq c, \label{eq:exactuW1p}
	\end{align}
	for $p=\infty$ if $d=2$, and $p=3$ if $d=3$.
	
	The same arguments applied on the boundary yield the corresponding estimates for $\psi_\ast$, with the advantage, that the dimension of the boundary is lower: 
	\begin{align}
		\| \psi_\ast^{(j)} \|_{L^\infty(\Ga_h)} + \| \widetilde \psi^k_\ast \|_{L^\infty(\Ga_h)} + \| \nabla_{\Ga_h} \widetilde \psi^k_\ast \|_{L^\infty(\Ga_h)}\leq c. \label{eq:exactpsiLinfty}
	\end{align}
	
	In the stability proof we then argue by a bootstrapping argument that, since the extrapolation only depends on past values, $\widetilde u$ is $L^\infty$-close to $\widetilde u_\ast$ and $\widetilde \psi$ is $L^\infty$-close to $\widetilde \psi_\ast$, hence all inputs are in a compact set. 
	
	\subsection{Results by Dahlquist and Nevanlinna $\&$ Odeh\protect\footnote{Up to some changes, the text of this preparatory subsection is taken from \cite[Section~6.3]{BullerjahnKovacs2024}.}}
	In order to derive energy estimates for BDF methods up to order $5$, we need the following important results from the $G$-stability theory of \cite[Theorem~3.3]{Dahlquist1978}, slightly extended to semi-inner products in \cite[Lemma~3.5]{ContriKovacsMassing2023}, and the multiplier technique of \cite[Section~2]{NevanlinnaOdeh1981}.
	%
	\begin{lem}
		\label{Lemma:Dahlquist} 
		Let $\delta (\zeta)= \sum_{j=0}^q \delta_j \zeta^j$ and $\mu(\zeta)=\sum_{j=0}^q \mu_j \zeta^j$ be polynomials of degree at most $q$, and at least one of them of degree $q$, that have no common divisor. Let $(\cdot , \cdot )$ be a semi-inner product on a complex Hilbert space $H$. 
		If
		\begin{align*}
			\Real \frac{\delta(\zeta)}{\mu(\zeta)} >0 , \qquad \text{ for all } \quad \zeta \in \mathbb{C}, \ |\zeta|<1,
		\end{align*}
		then there exists a symmetric positive definite matrix $G=(g_{ij}) \in \mathbb{R}^{q\times q}$, such that, for all $w_0,\dotsc,w_q \in H$, we have
		\begin{align*}
			\Real \Big( \sum_{i=0}^q \delta_i w_{q-i}, \sum_{j=0}^q \mu_j w_{q-j} \Big) \geq \sum_{i,j=1}^q g_{ij} ( w_i,w_j ) - \sum_{i,j=1}^q g_{ij} ( w_{i-1},w_{j-1} ) .
		\end{align*}
	\end{lem}
	
	The application of $G$-stability to the BDF schemes, for $q=1,\dotsc,5$, is ensured by the following result.
	\begin{lem}[{\cite[Section~2]{NevanlinnaOdeh1981}}]
		\label{Lemma:MultiplierTechnique}
		For $1\leq q\leq 5$, there exists $0 \leq \eta <1$ such that for $\delta (\zeta)= \sum_{j=0}^q \frac{1}{j} (1-\zeta)^j$,
		\begin{align*}
			\Real \frac{\delta(\zeta)}{1-\eta \zeta} >0 , \qquad \text{ for all } \quad \zeta \in \mathbb{C}, \ |\zeta|<1 .
		\end{align*}
		The classical values of $\eta$ from \cite[Table]{NevanlinnaOdeh1981} are found to be $\eta=0$, $0$, $0.0836$, $0.2878$, $0.8160$ for $q=1,\dotsc,5$ respectively.
	\end{lem}
	The exact multipliers were computed in \cite{AkrivisKatsoprinakis2015}, while multipliers for BDF-6 were derived in \cite{AkrivisChenYuZhou2021}.
	
	We thus introduce the $G$-semi-norm associated to the semi-inner product $(\cdot,\cdot)$ on a Hilbert space $H$: Given a collection of vectors $W^n=(w^n,\dotsc,w^{n-q+1}) \subset H^q$, we define
	\begin{align*}
		|W^n|_{G}^2:= \sum_{i,j=1}^q g_{ij} ( w^{n-i+1},w^{n-j+1} ),
	\end{align*}
	where $G$ is the symmetric positive definite matrix appearing in Lemma~\ref{Lemma:Dahlquist}. Then with the smallest and largest eigenvalues of $G$, denoted by $\lambda_0$ and $\lambda_1$, we have the inequalities:
	\begin{align}
		\lambda_0 |w^n|^2 \leq \lambda_0 \sum_{j=1}^q |w^{n-j+1}|^2 \leq |W^n|_G^2 \leq \lambda_1 \sum_{j=1}^q |w^{n-j+1}|^2, \label{eq:GnormEst}
	\end{align}
	where $|\cdot|$ is the semi-norm on $H$ induced by the semi-inner product $(\cdot,\cdot)$.
	Later on, an additional subscript, e.g. $|\cdot|_{G,a_h^{K,\alpha}}$, specifies which semi-inner product generates the $G$-weighted semi-norm.
	
	These results have previously been applied to the error analysis for the BDF time discretization of PDEs when testing the error equation with the error in \cite{LubichMansourVenkataraman2013}, and the discrete time derivative of the error in \cite{KovacsLiLubich2019}.

	\subsection{Stability result for $K \in [0,\infty)$ and $L\in (0,\infty)$} \label{sec:StabilityResult}
	
	In this section we give a detailed proof for the stability result of the discretized bulk--surface Cahn--Hilliard system with dynamic boundary conditions \eqref{eq:CHweakD}, for the majority of parameter choices $K \in [0,\infty)$ and $L\in (0,\infty)$, following Scheme~\ref{scheme:Stability-scheme}.
	
	\begin{prop} \label{prop:stability}
		Assume that the exact solution satisfies the regularity assumption \eqref{eq:regularityAss}, and $K\in [0,\infty)$ and $L\in(0,\infty)$. Consider the error equations of the numerical scheme \eqref{eq:erroreq}, for the linearly implicit BDF time discretization of order $1 \leq q \leq 5$. 
		Assume that, for step-sizes satisfying $\tau^q \leq C_0 h^\kappa$ (with arbitrary $C_0>0$), for a $\kappa > \frac{3}{2}$ the defects are bounded by 
		\begin{alignat*}{4}
			&\| (d^k_1,d^k_2) \|_{\ast,L,\beta} \leq ch^\kappa , \quad &&\| (d^k_3,d^k_4) \|_{\ast,K,\alpha} \leq ch^\kappa ,\\
			&\| (\partial^\tau d^i_3,\partial^\tau d^i_4) \|_{\ast,K,\alpha} \leq ch^\kappa, \quad &&\| (\partial^\tau_q d^k_3,\partial^\tau_q d^k_4) \|_{\ast,K,\alpha} \leq ch^\kappa,
		\end{alignat*}
		for $q\tau \leq k \tau \leq T$ and $1 \leq i \leq q-1$, where we recall that $\partial^\tau:=\partial_1^\tau$. 
		Further assume that the errors of the starting values satisfy
		\begin{align}
			\calE_I^{q-1}:= \sum_{i=0}^{q-1} \|(e_u^i,e_\psi^i)\|_{V_h^{K,\alpha}}^2  + \sum_{i=1}^{q-1} |\partial^\tau m_h\big((e_u^i,e_\psi^i),(\beta,1)\big)|^2\leq c h^{2\kappa} . \label{eq:DefIhq}
		\end{align}
		Then there exist $h_0 >0$ such that for all $h\leq h_0$ and $\tau >0$ with $\tau^q \leq C_0 h^\kappa$, the following stability estimate holds for \eqref{eq:erroreq}, for all $q \tau \leq t_n=n \tau \leq T$,
		\begin{align*}
			\| (e_u^n,e_\psi^n) \|_{V_h^{K,\alpha}}^2 + \tau \sum_{k=q}^n \|(e_\mu^k,e_\theta^k) \|_{V_h^{L,\beta}}^2 \leq C\Big(\calE_I^{q-1} + D_h^n\Big),
		\end{align*}
		where the constant $C>0$ is independent of $h$, $\tau$ and $n$, but depends exponentially on the final time $T$, and where 
		the term $D_h^n$ on the right-hand side collects the defects:
		\begin{align}
			D_h^n:=&\ \max_{q\leq k \leq n} \|(d_1^k,d_2^k) \|_{\ast,L,\beta}^2 +\|(d_3^n,d_4^n) \|_{\ast,K,\alpha}^2 + \|(d_3^q,d_4^q) \|_{\ast,K,\alpha}^2 +  \tau \sum_{k=q}^n \|(d_3^k,d_4^k) \|_{\ast,K,\alpha}^2 \nonumber \\
			&\ + \tau \sum_{k=q}^n \| (\partial^\tau_q d^k_3,\partial^\tau_q d^k_4) \|_{\ast,K,\alpha}^2 + \tau  \sum_{i=1}^{q-1} \|(\partial^\tau d^i_3,\partial^\tau d^i_4) \|_{\ast,K,\alpha}^2. \label{eq:DefDhn}
		\end{align}
	\end{prop}
	\begin{proof}
		Our strategy to obtain this stability bound is to use Scheme~\ref{scheme:Stability-scheme}.
		
		In the following $c$ and $C$ are generic positive constants, independent of $h$ and $\tau$, that may take different values on different occurrences. By $\rho >0$ we will denote a small number, independent of $h$ and $\tau$, used in Young's inequality, and hence we will often incorporate multiplicative constants into those yet unchosen factors. 
		
		We fix the constants $K\in [0,\infty)$, $L\in (0,\infty)$, $\alpha,\beta \in \R$, which satisfy $\alpha \beta |\Om| + |\Ga| \not= 0$ by assumption \hyperref[eq:AssumptionPara]{(A1)}, and $1\leq q \leq 5$. We denote by $0\leq \eta <1$ the multiplier factor as in Lemma~\ref{Lemma:MultiplierTechnique}, depending only on $q$. 
		
		\medskip
		
		\textit{(Preparations)} By our assumption \hyperref[eq:AssumptionPot]{(A2)}, $F_\Om'$, $F_\Ga'$, $F_\Om''$ and $F_\Ga''$ are \emph{locally Lipschitz continuous}. Since the extrapolation in the non-linear terms only depends on past values, we use a bootstrapping argument to show that $\widetilde u_h^n$ is $L^\infty$-close to $\widetilde u_\ast^n$ in $\Om_h$ and that $\widetilde \psi_h^n$ in $L^\infty$-close to $\widetilde \psi_\ast^n$ on $\Ga_h$, As indicated in Section~\ref{subsec:Linftyu}.
		
		Let $ t_{\max} \in (0,T]$ be the maximal time such that the following inequality holds for $k \tau \leq (n-1)\tau \leq t_{\max}$:
		\begin{align}
			\| e^{k}_u \|_{L^\infty(\Om_h)} + \| e_\psi^k\|_{L^\infty(\Ga_h)} \leq&\  h^{\frac{\kappa}{2}-\frac{d}{4}}. \label{eq:Linftybounde_u}
		\end{align}
		Note that for the starting values this bound is directly satisfied by the proposed bound, for $h$ sufficiently small. From now on we assume in the proof that $(n-1)\tau \leq t_{\max}$ and therefore these bounds are valid. At the end of the proof we then show that $t_{\max}=T$. It immediately follows, that for any $k = q+1,\dotsc,n$: 
		\begin{align*}
			\| \widetilde u_h^k-\widetilde u_\ast^k \|_{L^\infty(\Om_h)} + \|\widetilde \psi_h^k-\widetilde \psi_\ast^k \|_{L^\infty(\Ga_h)} \leq c h^{\frac{\kappa}{2}-\frac{d}{4}}.
		\end{align*}
		Together with \eqref{eq:exactuLinfty} and \eqref{eq:exactpsiLinfty}, we can conclude that there is a large enough compact set, which contains $\widetilde u_h^k$, $\widetilde u^k_\ast$, $\widetilde \psi_h^k$ and $\widetilde \psi_\ast^k$ (and even linear combinations), then on this set $F_\Om'$, $F_\Ga'$, $F_\Om''$ and $F_\Ga''$ are Lipschitz continuous by assumption.
		
		\textbf{Part A}: Let $q+1\leq k \leq n$, such that $(n-1)\tau \leq t_{\max}$. We get the energy estimate, by considering $\eqref{eq:erroreq1}^k$, i.e. equation \eqref{eq:erroreq1} considered at step $n=k$, tested with $(e_\mu^k-\eta e_\mu^{k-1},e_\theta^k-\eta e_\theta^{k-1}) \in V_h^{L,\beta}$ and subtracting $\eqref{eq:erroreq2}^k-\eta \eqref{eq:erroreq2}^{k-1}$ tested with $(\partial_q^\tau e_u^k, \partial_q^\tau e_\psi^k) \in V_h^{K,\alpha}$, this leads to the equality:
		\begin{align}
			&a_h^{K,\alpha}\big((\partial_q^\tau e_u^k, \partial_q^\tau e_\psi^k),(e_u^k-\eta e_u^{k-1},e_\psi^k-\eta e_\psi^{k-1})\big) + |(e_\mu^k,e_\theta^k)|_{a_h^{L,\beta}}^2 = \eta a_h^{L,\beta}\big((e_\mu^k,e_\theta^k),(e_\mu^{k-1},e_\theta^{k-1})\big) \nonumber \\
			&\ + m_h\big((d_1^k,d_2^k),(e_\mu^k-\eta e_\mu^{k-1},e_\theta^k-\eta e_\theta^{k-1})\big) - m_h\big((r_\Om^k-\eta r_\Om^{k-1},r_\Ga^k-\eta r_\Ga^{k-1}),(\partial_q^\tau e_u^k, \partial_q^\tau e_\psi^k)\big) \nonumber \\
			&\ - m_h\big((d_3^k-\eta d_3^{k-1},d_4^k-\eta d_4^{k-1}),(\partial_q^\tau e_u^k, \partial_q^\tau e_\psi^k)\big) \label{eq:testedEq}
		\end{align}
		
		\textbf{Part B}: We estimate the first term on the left-hand side of \eqref{eq:testedEq}, by applying Lemma~\ref{Lemma:Dahlquist} and Lemma~\ref{Lemma:MultiplierTechnique} to obtain, with the symmetric positive definite matrix $G$ and $E_u^k:=(e_u^k,\dots,e_u^{k-q+1})$, $E_\psi^k:=(e_\psi^k,\dots,e_\psi^{k-q+1})$:
		\begin{align}
			|(E_u^k,E_\psi^k)|_{G,a_h^{K,\alpha}}^2 - |(E_u^{k-1},E_\psi^{k-1})|_{G,a_h^{K,\alpha}}^2\leq \tau a_h^{K,\alpha}\big((\partial_q^\tau e_u^k, \partial_q^\tau e_\psi^k),(e_u^k-\eta e_u^{k-1},e_\psi^k-\eta e_\psi^{k-1})\big) .\label{eq:DahNOabs}
		\end{align}
		
		Inserting \eqref{eq:DahNOabs} into \eqref{eq:testedEq}, multiplying by $\tau$ and summing over $k=q+1,\dots,n$, we arrive at
		\begin{align}
			|(E_u^n,E_\psi^n)|_{G,a_h^{K,\alpha}}^2 + \tau \sum_{k=q+1}^n |(e_\mu^k,e_\theta^k)|_{a_h^{L,\beta}}^2 \leq&\ |(E_u^q,E_\psi^q)|_{G,a_h^{K,\alpha}}^2 + (I) + (II) + (III) + (IV)	\label{eq:EnergyEstUnbounded1}	
		\end{align}
		with the terms to be estimated in the following part:
		\begin{align*}
			(I) :=&\  \eta \tau \sum_{k=q+1}^n a_h^{L,\beta}\big((e_\mu^k,e_\theta^k),(e_\mu^{k-1},e_\theta^{k-1})\big), \\
			(II) :=&\ \tau \sum_{k=q+1}^n m_h\big((d_1^k,d_2^k),(e_\mu^k-\eta e_\mu^{k-1},e_\theta^k-\eta e_\theta^{k-1})\big), \\
			(III) :=&\ \tau \sum_{k=q+1}^n m_h\big((r_\Om^k-\eta r_\Om^{k-1},r_\Ga^k-\eta r_\Ga^{k-1}),(\partial_q^\tau e_u^k, \partial_q^\tau e_\psi^k)\big), \\
			(IV) :=&\ \tau \sum_{k=q+1}^n m_h\big((d_3^k-\eta d_3^{k-1},d_4^k-\eta d_4^{k-1}),(\partial_q^\tau e_u^k, \partial_q^\tau e_\psi^k)\big).
		\end{align*}
		\textbf{Part C}: We start by estimating the terms on the right-hand side of \eqref{eq:EnergyEstUnbounded1}. The term $(I)$ is estimated by Cauchy--Schwarz and Young's inequality, to obtain
		\begin{align}
			(I) \leq \eta \tau \sum_{k=q+1}^n |(e_\mu^k,e_\theta^k)|_{a_h^{L,\beta}}^2 + \frac{\eta}{2} \tau |(e_\mu^q,e_\theta^q)|_{a_h^{L,\beta}}^2. \label{eq:estCI}
		\end{align}
		
		The term $(II)$ is estimated by the definition of the dual norm \eqref{eq:DefectNorms} and Young's inequality, as
		\begin{align}
			(II) \leq \rho \tau \sum_{k=q+1}^n \| (e_\mu^k,e_\theta^k)\|_{V_h^{L,\beta}}^2 + c \tau \sum_{k=q+1}^n \| (d_1^k,d_2^k)\|_{\ast,L,\beta}^2 + c \tau \|(e_\mu^q,e_\theta^q)\|_{V_h^{L,\beta}}^2. \label{eq:estCII}
		\end{align}
		
		We consider the term $(III)$ and only the terms involving $(r_\Om^k,r_\Ga^k)$, since the terms in $(r_\Om^{k-1},r_\Ga^{k-1})$ can be treated in exactly the same way. Since $L\not= 0$ we can test \eqref{eq:erroreq1} with $(r_\Om^k,r_\Ga^k)$ and estimate as before:
		\begin{align}
			-m_h\big((r_\Om^k,r_\Ga^k),(\partial_q^\tau e_u^k, \partial_q^\tau e_\psi^k)\big) =&\ a_h^{L,\beta}\big((r_\Om^k,r_\Ga^k),(e_\mu^k,e_\theta^k)\big) - m_h\big((r_\Om^k,r_\Ga^k),(d_1^k,d_2^k)\big) \nonumber \\
			\leq&\ \rho |(e_\mu^k,e_\theta^k)|_{a_h^{L,\beta}}^2 + c \|(r_\Om^k,r_\Ga^k) \|_{V_h^{K,\alpha}}^2 + c \|(d_1^k,d_2^k)\|_{\ast,L,\beta}^2, \label{eq:testedL0}
		\end{align}
		where we have used in the last step the trace inequality, see, e.g., \cite[Section~5.5]{Evans1998}, and the fact that by definition, $\|\cdot \|_{V_h^{L,\beta}}:=a_h^\ast(\cdot,\cdot)$ does not depend on $L$ or $\beta$, see \eqref{eq:DefDiscNorms}.
		
		To bound the nonlinear terms $r_\Om^k$ and $r_\Ga^k$, we use the $H^1$-stability in order to remove the interpolant, see \cite[Theorem~5.9 \& Theorem~7.10]{ElliottRanner2021}, \eqref{eq:Linftybounde_u} and the \emph{preparations} and can therefore assume Lipschitz continuity of the nonlinear functions $F_\Om'$, $F_\Ga'$, $F_\Om''$ and $F_\Ga''$. We split this into the $L^2$-part, where we use the local Lipschitz continuity of $F_\Om'$ and $F_\Ga'$,
		\begin{align*}
			\|r_\Om^k\|_{H_h}^2 \leq  c \|F_\Om'(\widetilde u_h^k) - F_\Om'(\widetilde u_\ast^k) \|_{L^2(\Om_h)}^2 + c \|F_\Ga'(\widetilde \psi_h^k) - F_\Ga'(\widetilde \psi_\ast^k) \|_{L^2(\Ga_h)}^2 \leq c\| (\widetilde e_u^k,\widetilde e_\psi^k)\|_{H_h}^2,
		\end{align*}
		and the part involving the $L^2$-norm of the gradient, where we use product rule in space and Lipschitz continuity of $F_\Om''$ and $F_\Ga''$: 
		\begin{align*}
			\|(\nabla r_\Om^k,\nabla_{\Ga_h} r_\Ga^k)\|_{H_h}^2 \leq &\ c \|F_\Om''(\widetilde u_h^{k}) \nabla \widetilde u_h^{k} - F_\Om''(\widetilde u_\ast^{k})\nabla \widetilde u_\ast^{k} \|_{L^2(\Om_h)}^2 \nonumber \\
			&\ + c \|F_\Ga''(\widetilde \psi_h^{k}) \nabla_{\Ga_h} \widetilde \psi_h^{k} - F_\Ga''(\widetilde \psi_\ast^{k})\nabla_{\Ga_h} \widetilde \psi_\ast^{k} \|_{L^2(\Ga_h)}^2 \\
			=&\ c \| [F_\Om''(\widetilde u_h^{k})-F_\Om''(\widetilde u_\ast^{k})] \nabla \widetilde u_\ast^{k} + F_\Om''(\widetilde u_h^{k}) \nabla \widetilde e_u^{k}\|_{L^2(\Om_h)}^2 \nonumber \\
			&\ + c\| [F_\Ga''(\widetilde \psi_h^{k})-F_\Ga''(\widetilde \psi_\ast^{k})] \nabla_{\Ga_h} \widetilde \psi_\ast^{k} + F_\Ga''(\widetilde \psi_h^{k}) \nabla_{\Ga_h} \widetilde e_\psi^{k}\|_{L^2(\Ga_h)}^2 \\
			\leq&\ c \|\widetilde e_u^{k}\|_{L^q(\Om_h)}^2 \|\nabla \widetilde u_\ast^{k}\|_{L^p(\Om_h)}^2 +c \|\nabla \widetilde e_u^{k}\|_{L^2(\Om_h)}^2 \nonumber \\
			&\ + c \|\widetilde e_\psi^{k}\|_{L^2(\Ga_h)}^2 \|\nabla_{\Ga_h} \widetilde \psi_\ast^{k}\|_{L^\infty(\Ga_h)}^2  + c \|\nabla_{\Ga_h} \widetilde e_\psi^{k}\|_{L^2(\Ga_h)}^2\\
			\leq&\  c \|(\widetilde e_u^{k},\widetilde e_\psi^{k})\|_{V_h^{K,\alpha}}^2.
		\end{align*}
		For $d=2$ we estimate with $q=2$ and $p=\infty$ and use the $L^\infty$-bound for the Ritz map of the exact solutions \eqref{eq:exactuW1p} and \eqref{eq:exactpsiLinfty}, and for $d=3$ we choose $q=6$ and $p=3$ by H\"older's inequality and estimate further by Sobolev inequality (see e.g. \cite[Section~5.6.3]{Evans1998}), via the norm equivalence, and by the $L^3/L^\infty$-bound for the Ritz map of the exact solutions \eqref{eq:exactuW1p} and \eqref{eq:exactpsiLinfty}. 
		
		Altogether, we get for the term involving the non-linearity:
		\begin{align}
			(III) \leq&\ \rho \tau \sum_{k=q+1}^n |(e_\mu^k,e_\theta^k)|_{a_h^{L,\beta}}^2 +  c \tau \sum_{k=q+1}^{n-1} \| (e_u^k,e_\psi^k)\|_{V_h^{K,\alpha}}^2 + c \tau \| (e_u^q,e_\psi^q)\|_{V_h^{K,\alpha}}^2 \nonumber \\
			&\ + c \calE_I^{q-1} + c\tau \sum_{k=q+1}^n \|(d_1^k,d_2^k)\|_{\ast,L,\beta}^2. \label{eq:estCIII}
		\end{align}
		The term in $(e_u^q,e_\psi^q)$ will be estimated further later on, see \eqref{eq:qest}.
		
		The term $(IV)$ has to be treated differently, since we can not bound the $H^1$-norm of the defect, but also can not allow the term $(\partial_q^\tau e_u^k,\partial_q^\tau e_\psi^k)$ to appear on the right-hand side of \eqref{eq:EnergyEstUnbounded1}. The idea is to use a discrete product rule to transfer the discrete derivative on the error to a discrete derivative on the full term and on the defect. An adaptation of the analogous estimates in \cite[Equation~(10.20)]{KovacsLiLubich2019} and \cite[Equation~(6.38)]{BullerjahnKovacs2024} results in,
		\begin{align}
			(IV) \leq&\ \rho \sum_{k=n-q+1}^n \|(e_u^k,e_\psi^k)\|_{V_h^{K,\alpha}}^2 + \rho \tau \sum_{k=q+1}^n \|(e_u^k,e_\psi^k)\|_{V_h^{K,\alpha}}^2 + c \|(e_u^q,e_\psi^q)\|_{V_h^{K,\alpha}}^2 + c\calE_I^{q-1} \nonumber \\
			&\ + c \|(d_3^n,d_4^n)\|_{\ast,K,\alpha}^2 + c \|(d_3^q,d_4^q)\|_{\ast,K,\alpha}^2 + c\tau \sum_{k=q}^n \| (\partial_q^\tau d_3^k,\partial_q^\tau d_4^k)\|_{\ast,K,\alpha}^2 \nonumber \\
			&\  + c\tau \sum_{i=1}^{q-1} \| (\partial^\tau d_3^i,\partial_q^\tau d_4^i)\|_{\ast,K,\alpha}^2. \label{eq:estCIV}
		\end{align}
		
		Inserting now the estimates \eqref{eq:estCI}, \eqref{eq:estCII}, \eqref{eq:estCIII} and \eqref{eq:estCIV} into the energy estimate \eqref{eq:EnergyEstUnbounded1}, absorbing first the term in $\eta$, then the terms in $\rho$, by using the $G$-norm estimates \eqref{eq:GnormEst}, we arrive at
		\begin{align}
			\sum_{k=n-q+1}^n | (e_u^k,e_\psi^k) |_{a_h^{K,\alpha}}^2 + \tau \sum_{k=q+1}^n |(e_\mu^k,e_\theta^k)|_{a_h^{L,\beta}}^2 \leq&\ \rho \sum_{k=n-q+1}^n \|(e_u^k,e_\psi^k)\|_{V_h^{K,\alpha}}^2 \nonumber \\
			&\ + \rho \tau \sum_{k=q+1}^n \|(e_\mu^k,e_\theta^k)\|_{V_h^{L,\beta}}^2  \nonumber \\
			&\ + C \tau \sum_{k=q+1}^{n-1} \| (e_u^k,e_\psi^k) \|_{V_h^{K,\alpha}}^2 \nonumber \\
			&\ + C\big(\| (e_u^q,e_\psi^q) \|_{V_h^{K,\alpha}}^2 + \tau \|(e_\mu^q,e_\theta^q)\|_{V_h^{L,\beta}}^2\big) \nonumber \\
			&\ + C(\calE_I^{q-1}+ D_h^n). \label{eq:EnergyEstAlmost1}
		\end{align}
		
		\textbf{Part D}: Let $q\tau \leq n\tau \leq T$ and $(n-1)\tau \leq t_{\max}$, then the following Poincar\'e-type estimates hold,
		\begin{align}
			\| (e_u^n,e_\psi^n)\|_{V_h^{K,\alpha}} \leq&\ c | (e_u^n,e_\psi^n)|_{a_h^{K,\alpha}} + c \|(e_u^{q-1},e_\psi^{q-1})\|_{V_h^{K,\alpha}} + c \sum_{i=1}^{q-1} |\partial^\tau m_h((e_u^i,e_\psi^i),(\beta,1))| \nonumber \\
			&\ + c \max_{q\leq j \leq n}\|(d_1^j,d_2^j)\|_{\ast,L,\beta}, \label{eq:Poin1}\\
			\|(e_\mu^k,e_\theta^k) \|_{V_h^{L,\beta}} \leq&\ c | (e_\mu^n,e_\theta^n)|_{a_h^{L,\beta}} + c \| (\widetilde e_u^n,\widetilde e_\psi^n)\|_{V_h^{K,\alpha}}  + c \| (d_3^n,d_4^n)\|_{\ast,K,\alpha}.\label{eq:Poin2}
		\end{align}
		
		We start by showing the first inequality \eqref{eq:Poin1}, by using the geometric approximation results in Section~\ref{Subsection:GeomErrorEst}, to transfer $h$-independently between the discrete and continuous integrals, and then using the Poincar\'e--Wirtinger inequality for the continuous combined bulk and surface space $V^{K,\alpha}$, see Lemma~\ref{Lemma:BulkSurfPoincare}:
		\begin{align}
			\| (e_u^n,e_\psi^n)\|_{V_h^{K,\alpha}} \leq c  | (e_u^n,e_\psi^n)|_{a_h^{K,\alpha}} + c |\m_h^\beta(e_u^n,e_\psi^n)|, \label{eq:EstContPoin}
		\end{align}
		where we define the discrete combined mass $\m_h^\beta(e_u^n,e_\psi^n):=\beta |\Om_h| \langle e_u^n\rangle_{\Om_h} + |\Ga_h| \langle e_\psi^n \rangle_{\Ga_h}$. 
		
		Testing \eqref{eq:erroreq1} with $(\beta,1) \in V_h^{L,\beta}$ we obtain, with the notational convention $\m^n:=\m_h^\beta(e_u^n,e_\psi^n)$,
		\begin{align}
			\partial_q^\tau \m^n =m_h\big((d_1^n,d_2^n),(\beta,1)\big)=:\bfd^n. \label{eq:MassEq1}
		\end{align}
		Next we use, as in \cite[Eq.~(10.14)]{KovacsLiLubich2019}, the zero-stability of the BDF method to obtain the reformulation: 
		\begin{align}
			\partial^\tau \m^n=\sum_{k=q}^n \omega_{n-k} (\partial_q^\tau \m^k - \bfs^j), \quad \text{with} \quad \bfs^k=\begin{cases}
				0,\quad &k\geq2q-2 \\
				\sum_{i=1}^{q-1} \sigma_{k-i} \partial^\tau \m^i,\quad &\text{otherwise}
			\end{cases}, \label{eq:ReformPartial}
		\end{align}
		and where $|\omega_j|\leq c \rho^j$ for some $\rho<1$, and $\delta(\zeta)=(1-\zeta)\sum_{j=0}^{q-1}\sigma_j \zeta^j$, with $\sigma_j=0$ for $j\geq q$.
		
		The reformulation \eqref{eq:ReformPartial}, together with \eqref{eq:MassEq1}, yields the following estimate
		\begin{align*}
			\partial^\tau \m^n \leq&\ c \sum_{k=q}^n \rho^{n-k}  |\bfd^k| + c \sum_{i=1}^{q-1} |\partial^\tau \m^i| \nonumber \\
			\leq&\ c \frac{1}{1-\rho} \max_{q\leq k \leq n}|\bfd^k| + c \sum_{i=1}^{q-1} |\partial^\tau \m^i|.
		\end{align*}
		Then we obtain by telescopic summation:
		\begin{align*}
			|\m^n| \leq&\ |\m^{q-1}| + T \big(c \max_{q\leq k \leq n}|\bfd^k| + c \sum_{i=1}^{q-1} |\partial^\tau \m^i|\big) \nonumber \\
			\leq&\ c \|(e_u^{q-1},e_\psi^{q-1})\|_{V_h^{K,\alpha}} + c \sum_{i=1}^{q-1} |\partial^\tau \m^i| + c \max_{q\leq k \leq n}\|(d_1^k,d_2^k)\|_{\ast,L,\beta},
		\end{align*}
		which together with \eqref{eq:EstContPoin} results in the desired Poincar\'e-type estimate \eqref{eq:Poin1} for the phase fields.
		
		The second Poincar\'e type inequality \eqref{eq:Poin2} for the chemical potentials follows in a more direct way. We obtain, swapping the roles of $(K,\alpha)$ and $(L,\beta)$, in exactly the same way as for \eqref{eq:EstContPoin},
		\begin{align}
			\| (e_\mu^n,e_\theta^n) \|_{V_h^{L,\beta}} \leq c |(e_\mu^n,e_\theta^n)|_{a_h^{L,\beta}} + c |\m_h^\alpha(e_\mu^n,e_\theta^n)|, \label{eq:EstPoin2}
		\end{align} 
		
		Then we test \eqref{eq:erroreq2} with $(\alpha,1) \in V_h^{K,\alpha}$ and use local Lipschitz continuity and the \emph{preparations}, to obtain
		\begin{align}
			\m_h^\alpha(e_\mu^n,e_\theta^n) =&\ m_h\big((r_{\Om_h}^n,r_{\Ga_h}^n),(\alpha,1)\big) + m_h\big((d_3^n,d_4^n),(\alpha,1)\big) \nonumber \\
			\leq&\ c \|(\widetilde e_u^n,\widetilde e_\theta^n)\|_{V_h^{K,\alpha}} + c \|(d_3^n,d_4^n)\|_{\ast,K,\alpha} \label{eq:MassEst2},
		\end{align}
		and combining the estimates \eqref{eq:EstPoin2} and \eqref{eq:MassEst2} results in the Poincar\'e-type estimate \eqref{eq:Poin2} for the chemical potentials.
		
		\textbf{Part E}: We insert the Poincar\'e-type estimates \eqref{eq:Poin1} and \eqref{eq:Poin2} into inequality \eqref{eq:EnergyEstAlmost1}, and absorb terms to get:
		\begin{align*}
			\sum_{k=n-q+1}^n \| (e_u^k,e_\psi^k) \|_{V_h^{K,\alpha}}^2 + \tau \sum_{k=q+1}^n \|(e_\mu^k,e_\theta^k)\|_{V_h^{L,\beta}}^2 \leq&\  C \tau \sum_{k=q+1}^{n-1} \| (e_u^k,e_\psi^k) \|_{V_h^{K,\alpha}}^2 \nonumber \\
			&\ + C\big(\| (e_u^q,e_\psi^q) \|_{V_h^{K,\alpha}}^2 + \tau \|(e_\mu^q,e_\theta^q)\|_{V_h^{L,\beta}}^2\big) \nonumber \\
			&\ + C(\calE_I^{q-1}+ D_h^n). 
		\end{align*}
		The only terms which we need to further estimate are the terms of the time step $q$. Here we use the estimate \eqref{eq:qest} for the special case $n=q$. 
		
		Then we use a discrete Gr\"onwall inequality to arrive at the final estimate 
		\begin{equation}
			\label{eq:energy}
			\begin{aligned}
				\| (e_u^n,e_\psi^n) \|_{V_h^{K,\alpha}}^2 + \tau \sum_{k=q}^n \|(e_\mu^k,e_\theta^k)\|_{V_h^{L,\beta}}^2 \leq&\ C(\calE_I^{q-1}+ D_h^n).
			\end{aligned}
		\end{equation}
		
		\textit{(The special case $n=q$)} So far our analysis only covers the general case $n\geq q+1$. In the special case $n=q$ the analogous estimate
		\begin{align}
			\| (e_u^q,e_\psi^q) \|_{V_h^{K,\alpha}}^2 + \tau \|(e_\mu^q,e_\theta^q)\|_{V_h^{L,\beta}}^2 \leq&\ C(\calE_I^{q-1}+ D_h^q), \label{eq:qest}
		\end{align}	
		holds by a very similar proof. 
		
		The difference to the general case is that there is no equation for $n-1=q-1$. We solve this problem by adding the missing terms on both sides, such that we can use Dahlquist's $G$-stability result (Lemma~\ref{Lemma:Dahlquist}) and the multiplier technique of Nevanlinna and Odeh (Lemma~\ref{Lemma:MultiplierTechnique}), these terms involve only initial values, and therefore we estimated them directly.
		
		\textit{(Proving $t_{\max}=T$)}
		It is only left to show that $t_{\max}=T$ for sufficiently small $h \leq h_0$, see \eqref{eq:Linftybounde_u}. Suppose that $n$ is the maximal natural number such that $(n-1) \tau \leq t_{\max}$ holds and $n \tau \leq T$. We then have, by the above inequality \eqref{eq:energy} and the assumed bounds on the defects and initial values, that $\| (e^n_u, e_\psi^n)\|_{V_h^{K,\alpha}} \leq Ch^\kappa$. 
		
		Then by the inverse estimate (see \cite[Theorem~4.5.11]{BrennerScott2008}), we have that
		\begin{align*}
			\| e^n_u \|_{L^\infty(\Om_h)} \leq&\  ch^{-\frac{d}{2}} \| e^n_u \|_{L^2(\Om_h)} \leq ch^{-\frac{d}{2}} \| (e^n_u, e_\psi^n)\|_{V_h^{K,\alpha}} \leq C h^{\kappa-\frac{d}{2}} \leq h^{\frac{\kappa}{2}-\frac{d}{4}}, \\
			\| e^n_\psi \|_{L^\infty(\Ga_h)} \leq&\  ch^{-\frac{d-1}{2}} \| e^n_\psi \|_{L^2(\Ga_h)} \leq ch^{-\frac{d-1}{2}} \| (e^n_u, e_\psi^n)\|_{V_h^{K,\alpha}} \leq C h^{\kappa-\frac{d-1}{2}} \leq h^{\frac{\kappa}{2}-\frac{d}{4}},
		\end{align*}
		for sufficiently small $h$. Therefore $n\tau \leq t_{\max}$, which contradicts the maximality of $t_{\max}$, hence $t_{\max}=T$.	\hfill	 
	\end{proof}
	
	\subsection{Stability results for special cases of the parameters $K$ and $L$}
	
	The stability result Proposition~\ref{prop:stability} for the case $K\in [0,\infty)$ and $L\in (0,\infty)$ is extended to the cases $K \in (0,\infty)$ and $L=\infty$, as well as $K=L=0$, by minimal changes in the above presented proof.
	
	\subsubsection{Stability result in the case $K \in (0,\infty)$, $L=\infty$}
	
	\begin{prop} \label{prop:stab2}
		The stability result given in Proposition~\ref{prop:stability} also holds in the case $K \in (0,\infty)$, $L=\infty$.
	\end{prop}
	\begin{proof}
		The proof is identical to the proof of Proposition~\ref{prop:stability}, with the exception of \textbf{Part D}. 
		
		Here the individual mass conservation on bulk and surface for the phase field $(u,\psi)$, 
		\begin{align*}
			\int_\Om u(t) = \int_\Om u_0 \text{ and } \int_\Ga \psi(t) = \int_\Ga \psi_0,
		\end{align*}
		induced by the vanishing mass flux $\partial_\nu \mu=0$ (by the choice $L=\infty$), leads along the same lines as in the original \textbf{Part D} to an individual \emph{almost mass conservation} of the error in $u$ and $\psi$. Just instead of testing with $(\beta,1)$ we now test with $(1,0)$ and $(0,1)$. This individual \emph{almost mass conservation} then allows us to use the classical Poincar\'e--Wirtinger inequalities on bulk and surface separately, see, e.g., \cite[Section~5.8.1]{Evans1998}, \cite[Theorem~2.12]{DziukElliott2013}, which then leads to the same estimate \eqref{eq:Poin1}. The analogous estimate to \eqref{eq:Poin2} follows along the same lines, note that here we need $K>0$ to be able to test with $(1,0)$ and $(0,1)$.
	\end{proof}
	
	\subsubsection{Stability result in the case $K =0$, $L=0$}
	
	\begin{prop} \label{prop:stab3}
		The stability result given in Proposition~\ref{prop:stability} also holds in the case $K =0$, $L=0$.
	\end{prop}
	\begin{proof}
		The proof is identical to the proof of Proposition~\ref{prop:stability}, with the exception of \eqref{eq:testedL0}, where we are only allowed to test \eqref{eq:erroreq1} with $(r_\Om^n,r_\Ga^n) \in V_h^{L,\beta}$, for $L=0$, if $\displaystyle{ \beta r_\Ga^n = r_\Om^n|_{\Ga_h}}$. This is assured by the additional condition in \hyperref[eq:AssumptionPara]{(A1)} for the parameter choice $K=L=0$.  
	\end{proof}

	\section{Consistency \label{section:Consistency}}
	The consistency analysis relies on a combination of the error analysis for the spatial discretization, which results from an interplay of the geometric errors and the Ritz map error estimates, similar to \cite[Section~6]{HarderKovacs2022}, and the temporal error estimates, using the techniques developed in \cite[Section~7]{BullerjahnKovacs2024}. In the following section we prove optimal-order bounds for the fully discrete defect terms.

	\subsection{Defect bounds}
	Using the geometric approximation result from Section~\ref{Subsection:GeomErrorEst} and temporal error estimates from Section~\ref{sec:TempErrorEst}, we arrive at the optimal-order bounds for the different defect terms which appear in the stability analysis.
	\begin{prop} \label{prop:consistency}
		Let $(u,\psi,\mu,\theta)$ solve \eqref{eq:CHweak} and satisfy the regularity assumptions \eqref{eq:regularityAss}. Then we have, for $q \leq k \leq n$ and $1\leq i \leq q-1$,
		\begin{alignat*}{3}
			\| (d^k_1,d_2^k) \|_{\ast,L,\beta} + \| (d^k_3,d_4^k) \|_{\ast,K,\alpha} + \| (\partial^\tau_q d^k_3,\partial_q^\tau d_4^k) \|_{\ast,K,\alpha} \leq&\  C (h^2 + \tau^q), \\
			\| (\partial^\tau d^i_3,\partial^\tau d_4^i) \|_{\ast,K,\alpha} \leq&\   C (h^2 + \tau^q) .
		\end{alignat*}
		With these estimates we obtain the bound on the full defect term, defined in \eqref{eq:DefDhn},
		\begin{align*}
			D_h^n \leq C (h^2 + \tau^{q})^2.
		\end{align*}
	\end{prop}
	\begin{proof}
		The proof is an adaptation of the optimal-order defect bounds in \cite{HarderKovacs2022,BullerjahnKovacs2024} to the present case. 
		
		First, we obtain optimal bounds for the semi-discrete defects and their first and second time derivative, following the lines of \cite[Section~7]{HarderKovacs2022}: We add the exact weak formulation \eqref{eq:CHweak} to the defect as a zero, exchange the bilinear forms involving gradients to $L^2$-bilinear forms by the definition of the Ritz map and use the geometric error estimates in Lemma~\ref{Lemma:GeomAppError} and Ritz map error estimates in Lemma~\ref{Lemma:Ritzmaperror} to obtain bounds of optimal order. The bounds for the terms involving the nonlinear potentials additionally need the $L^\infty$-bounds in Section~\ref{subsec:Linftyu} and the assumed local Lipschitz continuity \hyperref[eq:AssumptionPot]{(A2)}, as well as standard estimates for the finite element interpolation, see \cite[Theorem~5.9 \& Theorem~7.10]{ElliottRanner2021}.
		
		Then we use Lemma~\ref{Lemma:PeanoKernel} to conclude the optimal-order defect bounds for $(d_1^k,d_2^k)$ and $(d_3^k,d_4^k)$. The estimate on the discrete derivative of the defect term $(d_3^k,d_4^k)$ is derived following closely the lines of \cite[Section~7.3]{BullerjahnKovacs2024}, by introducing an extension of the definition of the extrapolation to the interval $[0,q\tau]$ by a quintic Hermite polynomial and using Lemma~\ref{Lemma:PeanoKernel}. As an immediate consequence of the discrete defect term $(d_3^k,d_4^k)$ only containing the extrapolation and not the discrete derivative, we can allow the exact solution to be one derivative less regular than in \cite{BullerjahnKovacs2024}. The estimates in $(\partial^\tau d^i_3,\partial^\tau d_4^i)$ are an immediate consequence of the semi-discrete defect bounds.
	\end{proof}

	\section{Optimal order fully discrete error estimates}
	\label{section:MainProof}
	We are now ready to combine the stability and consistency analysis of the last two sections to establish the optimal-order fully discrete error estimates of the numerical scheme \eqref{eq:CHweakD} for the bulk--surface Cahn--Hilliard system with dynamic boundary conditions, which we formulated in Theorem~\ref{thm:OptimalOrderErrorEst}.
	
	\begin{proof}[Proof of Theorem~\ref{thm:OptimalOrderErrorEst}]
		We use the standard method of decomposing the error into the projection error and the error of the numerical scheme and the projection of the exact solution at time $t_n$, see, e.g., \cite{Thomee2006}:
		\begin{equation} \label{eq:MainProof1}
			\begin{aligned}
				(u^n_h)^\ell - u(t_n)=&\  (e^n_u)^\ell + R_h^\Om u(t_n) - u(t_n).
			\end{aligned}
		\end{equation}
		
		For the second term we use the Ritz map error estimates of Lemma~\ref{Lemma:Ritzmaperror} and the regularity assumptions \eqref{eq:regularityAss} to estimate
		\begin{equation} \label{eq:MainProof2}
			\begin{aligned}
				\| R_h^\Om u - u \|_{L^2(\Om)}  + h \| R_h^\Om u - u \|_{H^1(\Om)} \leq&\ ch^2,
			\end{aligned}
		\end{equation}
		
		For the first term we use a combination of the stability part, Propositions~\ref{prop:stability}, \ref{prop:stab2}, \ref{prop:stab3}, with the bounds for the defects in Proposition~\ref{prop:consistency} and the mild step size restriction $\tau^q\leq C_0 h^2$, to arrive at
		\begin{align}
			\Big(\| (e^n_u)^\ell\|_{H^1(\Om)}^2 +\| (e^n_\psi)^\ell\|_{H^1(\Ga)}^2 + \tau \sum_{k=q}^n \| (e^k_\mu)^\ell \|_{H^1(\Om)}^2 + \| (e^k_\theta)^\ell \|_{H^1(\Ga)}^2 \Big)^{\frac{1}{2}} \leq C (h^2+\tau^q), \label{eq:MainProof3}
		\end{align}
		using again the norm equivalence, see Remark~\ref{rem:normequivalence}. Combining all above estimates \eqref{eq:MainProof1}, \eqref{eq:MainProof2} and \eqref{eq:MainProof3} gives the stated result for the bulk phase field $u$. 
		
		The estimates for $\psi$, $\mu$ and $\theta$ follow in exactly the same way.
	\end{proof}

	\begin{rem}
		We strongly believe that Theorem~\ref{thm:OptimalOrderErrorEst} can be extended to higher-order isoparametric bulk--surface finite elements, by straightforward, yet technically non-trivial arguments. Most of the changes are restricted to Section~\ref{section:Consistency}, and sporadically elsewhere, whenever spatial approximation results are used. The analogous Ritz map and geometric approximation results can be found in \cite{Demlow2009}, \cite{Kovacs2018} and \cite{ElliottRanner2021}.
	\end{rem}

	\section{Generalization and the Cahn--Hilliard equation on evolving surfaces \label{section:Generalization}}
	
	As indicated in the previous analysis we expect that the general framework of this optimal-order convergence result, and especially the proof of the stability result using Scheme~\ref{scheme:Stability-scheme}, translates easily to other problems satisfying an energy estimate and \emph{almost} mass conservation properties. 
	
	Such equations are for instance the Allen--Cahn and Cahn--Hilliard equations in stationary domains with classical boundary conditions, see \cite{AllenCahn1979}, respectively \cite{CahnHilliard1958}, or the same PDEs on evolving surfaces, see \cite{OlshanskiiXuYushutin2021}, respectively \cite{ElliottRanner2015,CaetanoElliott2021}. 
	
	To demonstrate the generality of the presented approach, let us show how these arguments transfer over to the Cahn--Hilliard equation on evolving surfaces, see, e.g. \cite{ElliottRanner2015,CaetanoElliott2021}:
	\begin{subequations} \label{eq:CHevolSurf}
		\begin{gather}
			\matdev u + u(\nabla_\Ga \cdot V) - \Delta_\Ga w =0, \\
			w + \epsilon \Delta_\Ga u = \frac1\epsilon F_\Ga'(u),
		\end{gather}
	\end{subequations}
	subject to the suitable initial condition $u(0)=u_0$, and where $\matdev$ denotes the material derivative and $V$ the velocity field of the surface, for more details see, e.g., \cite{ElliottRanner2015}.
	
	We base our analysis on the evolving surface finite element method developed in \cite{DziukElliott2013}. This spatial discretization for the evolving surface Cahn--Hilliard equation was first analyzed in \cite{ElliottRanner2015}, and slightly different semi-discrete error estimates were shown in \cite{BeschleKovacs2022}. In \cite{ElliottSales2024B,ElliottSales2024} a backward Euler time discretization was used to derive optimal-order fully discrete error estimates, for a wide range of non-linear potentials.
	
	Following standard notation in the literature we denote by $\Ga_h(t)$ the triangulated surface and by $S_h(t)$ the surface finite element space at time $t \in [0,\infty)$. 
	
	Using this spatial discretization we define the bilinear forms for all $t\in [0,\infty)$,
	\begin{align*}
		&m_{\Ga_h}(t;\cdot,\cdot) \colon S_h(t) \times S_h(t) \to \mathbb{R} ,\qquad &&m_{\Ga_h}(t;\phi_h,\psi_h) := \int_{\Ga_h(t)} \phi_h \psi_h , \\
		&a_{\Ga_h}(t;\cdot,\cdot) \colon S_h(t) \times S_h(t) \to \mathbb{R} ,\qquad
		&& a_{\Ga_h}(t;\phi_h,\psi_h) :=  \int_{\Ga_h(t)} \nabla_{\Ga_h(t)} \phi_h \cdot \nabla_{\Ga_h(t)} \psi_h  ,
	\end{align*}
	where we usually suppress the dependence on a timestep $t_n$ in the superscript $n$, when considering the full discretization. 
	
	We combine this with the linearly implicit backward difference method of order $q=1,...,5$ in time to arrive at the numerical scheme: Find for $h,\tau >0$ and $q \tau \leq t_n:=n\tau \leq T$, functions $u_h^n,w_h^n \in S_h^n:=S_h(t_n)$, such that
	\begin{subequations}\label{eq:SurfNumScheme}
		\begin{align}
			\partial_q^\tau \Big( m_{\Ga_h}^n(u_h^n,\phi_h^n)\Big) + a_{\Ga_h}^n(w_h^n,\phi_h^n) =&\ 0, \\
			m_{\Ga_h}^n(w_h^n,\psi_h^n) - \epsilon a_{\Ga_h}^n(u_h^n,\psi_h^n) =&\ \frac1\epsilon m_{\Ga_h}^n(\calI_{\Ga_h^n} F_\Ga'(\widetilde u_h^n),\psi_h^n),
		\end{align}
	\end{subequations}
	for all $\phi_h^n, \psi_h^n \in S_h^n$, and given the starting values $u_h^{0} \in S_h^{0},\cdots,u_h^{q-1} \in S_h^{q-1}$. Here we define 
	\begin{align*}
		\partial_q^\tau \Big( m_{\Ga_h}^n(u_h^n,\phi_h^n)\Big):= \frac1\tau \sum_{i=0}^q \delta_i m_{\Ga_h}^{n-i}(u_h^{n-i},\underline{\phi_h^{n}}), \qquad 	\widetilde u_h^n:= \sum_{i=1}^q \gamma_i \overline{u_h^{n-i}},
	\end{align*}
	with definition \eqref{eq:BdfCoeff}, and where $\underline{\phi_h^{n}} \in S_h^{n-i}$ denotes the function in $S_h^{n-i}$ with the same nodal values as $\phi_h^n$, and analogously $\overline{u_h^{n-i}} \in S_h^n$ denotes the function in $S_h^n$ with the same nodal values as $u_h^{n-i} \in S_h^{n-i}$. 
	
	For a sufficiently regular solution pair $(u,w)$, this leads to the error terms $e_u^n=u_h^n-\widetilde R_h^\Ga u(t_n)$ and $e_w^n=w_h^n - \widetilde R_h^\Ga w(t_n)$, for the Ritz map $\widetilde R_h^\Ga$ defined in \cite[Definition~2.18]{ElliottSales2024}. Finally, we obtain, with $\epsilon=1$, the error equations
	\begin{subequations} \label{eq:SurfError}
		\begin{align}
			\partial_q^\tau \Big( m_{\Ga_h}^n(e_u^n,\phi_h^n)\Big) + a_{\Ga_h}^n(e_w^n,\phi_h^n) =&\ - m_{\Ga_h}^n(d_1^n,\phi_h^n), \label{eq:SurfError1}\\
			m_{\Ga_h}^n(e_w^n,\psi_h^n) - \epsilon a_{\Ga_h}^n(e_u^n,\psi_h^n) =&\ m_{\Ga_h}^n(r_h^n,\psi_h^n) - m_{\Ga_h}^n(d_2^n,\psi_h^n),\label{eq:SurfError2}
		\end{align}
	\end{subequations}
	where we denote the defect terms as $d_1^n$ and $d_2^n$, and collect the terms involving the nonlinear potential in $r_h^n = \calI_{\Ga_h^n} F_\Ga'(\widetilde u_h^n) - \calI_{\Ga_h^n} F_\Ga'(\widetilde u_\ast^n)$, where $\widetilde u_\ast^n:=\widetilde{ \widetilde R_h^\Ga u(t_n)}$.

	\subsection{Stability estimate}
	
	The major difference in the stability analysis stems from the dependence on time of the bilinear forms, which can be handled by the following result from \cite[Lemma~3.4]{ElliottSales2024}, originating from \cite[Section~3]{DziukElliott2012}:
	
	\begin{lem} \label{lem:SurfDifference}
		For $\zeta_h^n,\eta_h^n \in S_h^n$ and sufficiently small $\tau$ we have 
		\begin{align*}
			|\mh^n(\zeta_h^n,\eta_h^n) - \mh^{n-1}(\underline{\zeta_h^n},\underline{\eta_h^n})| \leq&\ c \tau \| \zeta_h^n \|_{L^2(\Ga_h^n)} \| \eta_h^n \|_{L^2(\Ga_h^n)}, \nonumber \\
			|\ah^n(\zeta_h^n,\eta_h^n) - \ah^{n-1}(\underline{\zeta_h^n},\underline{\eta_h^n})| \leq&\ c \tau \| \nabla_{\Ga_h} \zeta_h^n \|_{L^2(\Ga_h^n)} \| \nabla_{\Ga_h} \eta_h^n \|_{L^2(\Ga_h^n)},
		\end{align*}
		where $c>0$ denotes a constant independent of $h$ and $\tau$.
	\end{lem}
	Note that this result easily extends to the case of up to $q$ time-steps. Additionally, this provides a norm equivalence for $L^2$- and $H^1$-norms on different surfaces, see also \cite[Lemma~3.6]{DziukElliott2012}, \cite[Corollary~3.5]{ElliottSales2024}.
	
	Using Scheme~\ref{scheme:Stability-scheme}, we obtain a similar stability estimate as in Proposition~\ref{prop:stability}.
	
	\begin{prop} \label{prop:Surfstability}
		Assume that the exact solution satisfies the regularity assumption \eqref{eq:SurfregularityAss}. Consider the error equations of the numerical scheme \eqref{eq:SurfError}, for the linearly implicit BDF time discretization of order $1 \leq q \leq 5$. 
		Assume that, for step-sizes satisfying $\tau^q \leq C_0 h^\kappa$ (with arbitrary $C_0>0$), for a $\kappa > 1$ the defects are bounded by 
		\begin{align}
			D_h^n:=&\ \max_{q\leq k \leq n} \|d_1^k \|_{\ast,k}^2 + \max_{q\leq k \leq n} \|d_2^k \|_{\ast,k}^2 + \tau \sum_{k=q}^n \| \partial^\tau_q d^k_2\|_{\ast,k}^2 + \tau  \sum_{i=1}^{q-1} \|\partial^\tau d^i_2\|_{\ast,k}^2 \leq c h^{2\kappa}, \label{eq:SurfDefDhn}
		\end{align}
		where we denote by $\| \cdot \|_{\ast,k}$ the dual norm to $H^1(\Ga_h^k)$, and recall that $\partial^\tau:=\partial_1^\tau$. Further assume that the errors of the starting values satisfy
		\begin{align}
			\calE_I^{q-1}:= \sum_{i=0}^{q-1} \| e_u^i\|_{H^1(\Ga_h^i)}^2  + \sum_{i=1}^{q-1} |\partial^\tau m_{\Ga_h}^n\big(e_u^i,1\big)|^2\leq c h^{2\kappa} . \label{eq:SurfDefIhq}
		\end{align}
		Then there exist $h_0 >0$ such that for all $h\leq h_0$ and $\tau >0$ with $\tau^q \leq C_0 h^\kappa$, the following stability estimate holds for \eqref{eq:SurfError}, for all $q \tau \leq t_n=n \tau \leq T$,
		\begin{align*}
			\| e_u^n\|_{H^1(\Ga_h^n)}^2 + \tau \sum_{k=q}^n \|e_w^k\|_{H^1(\Ga_h^k)}^2 \leq C\Big(\calE_I^{q-1} + D_h^n\Big),
		\end{align*}
		where the constant $C>0$ is independent of $h$, $\tau$ and $n$, but depends exponentially on the final time $T$.
	\end{prop}
	\begin{proof}
		
		\textit{(Preparations)} We show that there is a large enough compact set, which contains $\widetilde u_h^k$ and $\widetilde u_\ast^k$, such that we can assume that $F_\Ga'$ and $F_\Ga''$ are Lipschitz continuous on this set. Let $t_{\max} \in (0,T]$ be the maximal time such that the following inequality holds for $k \tau \leq (n-1)\tau \leq t_{\max}$:
		\begin{align}
			\| e_u^k \|_{L^\infty(\Ga_h^k)} \leq h^{\frac\kappa2 - \frac12}. \label{eq:Surf Bootstrap}
		\end{align}
		Note that for the starting values this bound is directly satisfied for $h$ sufficiently small. Together with an $L^\infty$-bound for the Ritz map, analogous to Section~\ref{sec:Ritzmap}, we obtain the proposed set and we can assume Lipschitz continuity throughout the proof for $k \tau \leq t_{\max}$. It remains to show that $t_{\max} = T$, which we obtain by a bootstrapping argument at the end of the proof.
		
		\textbf{Part A}: Let $q+1 \leq k \leq n$ and test equation $\eqref{eq:SurfError1}^k$ with $(e_w^k - \eta \overline{e_w^{k-1}})$, subtract equation $\eqref{eq:SurfError2}^k$ tested with $\partial_q^\tau e_u^k$ and add equation $\eta \eqref{eq:SurfError2}^{k-1}$ tested with $\underline{\partial_q^\tau e_u^k}$, to obtain
		\begin{align*}
			\ah^k(e_u^k - \eta \overline{e_u^{k-1}},\partial_q^\tau e_u^k) + \ah^k(e_w^k,e_w^k) = \sum_{i=1}^6 R_i^k , 
		\end{align*}
		where 
		\begin{align*}
			R_1^k :=&\ \eta \ah^k(e_w^k,\overline{e_w^{k-1}}), \\
			R_2^k :=&\ - \mh^k(d_1^k,e_w^k - \eta \overline{e_w^{k-1}}), \\
			R_3^k :=&\ \mh^k(r_h^k,\partial_q^\tau e_u^k) - \eta \mh^{k-1}(r_h^{k-1}, \underline{\partial_q^\tau e_u^k}), \\
			R_4^k :=&\ - \Big[\mh^k(d_2^k,\partial_q^\tau e_u^k) - \eta \mh^{k-1}(d_2^{k-1},\underline{\partial_q^\tau e_u^k})\Big], \\
			R_5^k :=&\ \eta \Big[\ah^k(\overline{e_u^{k-1}},\partial_q^\tau e_u^k) - \ah^{k-1}(e_u^{k-1},\underline{\partial_q^\tau e_u^k})\Big], \nonumber \\
			R_6^k :=&\ \Big[\mh^k(e_w^k,\partial_q^\tau e_u^k) - \eta \mh^{k-1}(e_w^{k-1}, \underline{\partial_q^\tau e_u^k})\Big] - \frac1\tau \sum_{i=0}^q \delta_i \mh^{k-i}(e_u^{k-i},e_w^k - \eta \overline{e_w^{k-1}}).
		\end{align*}
		From now on we simplify the notation, e.g. $\overline{e_u^{k-1}}$ to $e_u^{k-1}$, as it is always clear from the context on which domain the finite element function is considered.

		\textbf{Part B}: We estimate the first term on the left-hand side of \eqref{eq:Surf Eq tested}, by applying Lemma~\ref{Lemma:Dahlquist} and Lemma~\ref{Lemma:MultiplierTechnique} to obtain, with the symmetric positive definite matrix $G$ and $E_u^k:=(e_u^k,\dots,e_u^{k-q+1})$:
		\begin{align*}
			|E_u^k|_{G,\ah^k}^2 - |E_u^{k-1}|_{G,\ah^{k-1}}^2 \leq \tau \ah^k(e_u^k - \eta e_u^{k-1},\partial_q^\tau e_u^k) + \Big(|E_u^{k-1}|_{G,\ah^{k}}^2 - |E_u^{k-1}|_{G,\ah^{k-1}}^2\Big)
		\end{align*}
		Summing over $k=q+1,\dots,n$ and multiplying by $\tau$, we obtain
		\begin{align}
			|E_u^n|_{G,\ah^n}^2 + \tau \sum_{k=q+1}^n |e_w^k|_{\ah^k}^2 \leq&\  |E_u^q|_{G,\ah^q}^2 + \tau \sum_{k=q+1}^n L_1^k + \tau \sum_{k=q+1}^n \sum_{i=1}^6 R_i^k, \label{eq:Surf Eq tested}
		\end{align}
		where we have from the left-hand side,
		\begin{align*}
			L_1^k:=&\ \frac1\tau \Big(|E_u^{k-1}|_{G,\ah^{k}}^2 - |E_u^{k-1}|_{G,\ah^{k-1}}^2\Big).
		\end{align*}
		
		\textbf{Part C}: We start by estimating the terms on the right-hand side. The estimates are usually split into a part containing the surface discrepancy depending on time, which is estimated by Lemma~\ref{lem:SurfDifference}, and a remaining part, which is estimated analogously to \textbf{Part C} in Proposition~\ref{prop:stability}. In the following estimates, we collect terms in the defects, initial data and terms which can be treated by Gr\"onwall's inequality and the case $n=q$ in
		\begin{equation*}
			\kappa_h^n:= \tau \sum_{k=q+1}^n \| e_u^k \|_{H^1(\Ga_h^k)}^2 + \| e_u^q \|_{H^1(\Ga_h^q)}^2 + \tau  \| e_w^q \|_{H^1(\Ga_h^q)}^2 + \calE_I^{q-1} + D_h^n.
		\end{equation*}  
		
		We start by estimating the term in $L_1$ by Lemma~\ref{lem:SurfDifference} and Young's inequality, as
		\begin{align}
			\tau \sum_{k=q+1}^n L_1^k \leq&\ \sum_{k=q+1}^n \sum_{i,j=0}^{q-1} |g_{ij}| |\ah^k(e_u^{k-i-1},e_u^{k-j-1}) - \ah^{k-1}(e_u^{k-i-1},e_u^{k-j-1})| \nonumber \\
			\leq&\ c \tau \sum_{k=q+1}^n \sum_{i,j=0}^{q-1} | e_u^{k-i-1} |_{\ah^{k-i-1}} | e_u^{k-j-1} |_{\ah^{k-j-1}} \nonumber \\
			\leq&\ c  \kappa_h^n. \label{eq:Surf C L1}
		\end{align}
		We note that here and in the following, we use in almost every step the norm equivalence of norms on different surfaces, see \cite[Lemma~3.6]{DziukElliott2012}, \cite[Corollary~3.5]{ElliottSales2024}.
		
		The term in $R_1$ is estimated again by Lemma~\ref{lem:SurfDifference}, as
		\begin{align}
			\tau \sum_{k=q+1}^n R_1^k \leq&\ \eta \tau \sum_{k=q+1}^n | e_w^{k} |_{\ah^{k}} \big(| e_u^{k-1} |_{\ah^{k-1}} + c \tau | e_u^{k-1} |_{\ah^{k-1}}\big) \nonumber \\
			\leq&\ \eta \tau \sum_{k=q+1}^n | e_w^{k} |_{\ah^{k}}^2 + \delta \tau \sum_{k=q+1}^n | e_w^{k} |_{\ah^{k}}^2 +  c  \kappa_h^n, \label{eq:Surf C R1}
		\end{align}
		for $\tau$ sufficiently small.
		
		The term in $R_2$ is estimated by Cauchy-Schwartz and Young's inequality and the definition of the dual norm, as
		\begin{align}
			\tau \sum_{k=q+1}^n R_2^k \leq&\ \rho \tau \sum_{k=q+1}^n \| e_w^{k} \|_{H^1(\Ga_h^k)}^2 +  c  \kappa_h^n.\label{eq:Surf C R2}
		\end{align}
		
		For the term in $R_3$ we aim to use \eqref{eq:SurfError1} to eliminate the derivative, as is done in \textbf{Part C} of Proposition~\ref{prop:stability}. In order to do this, we again use Lemma~\ref{lem:SurfDifference} and the \emph{preparations} to estimate:
		\begin{align}
			\tau \sum_{k=q+1}^n R_2^k \leq&\ \tau \sum_{k=q+1}^n \mh^k(r_h^k -\eta r_h^{k-1},\partial_q^\tau e_u^k) \nonumber \\
			&\ + \eta \tau \sum_{k=q+1}^n \Big[\mh^k(r_h^{k-1},\partial_q^\tau e_u^k) - \mh^{k-1}(r_h^{k-1},\partial_q^\tau e_u^k)\Big] \nonumber \\
			\leq&\ \tau \sum_{k=q+1}^n \mh^k(r_h^k -\eta r_h^{k-1},\partial_q^\tau e_u^k) \nonumber \\
			&\ + c \tau \sum_{k=q+1}^n \| r_h^{k-1} \|_{L^2(\Ga_h^{k-1})} \tau \|\partial_q^\tau e_u^{k} \|_{L^2(\Ga_h^k)} \nonumber \\
			\leq&\ \tau \sum_{k=q+1}^n \mh^k(r_h^k -\eta r_h^{k-1},\partial_q^\tau e_u^k) + c \kappa_h^n , \label{eq:Surf Helper derivative}
		\end{align}
		then we add in each summand of the first term $\eqref{eq:SurfError1}^k$ tested with $r_h^k - \eta r_h^{k-1}$, to obtain
		\begin{align*}
			\tau \sum_{k=q+1}^n R_2^k \leq&\ \tau \sum_{k=q+1}^n \Big[\mh^k(r_h^k -\eta r_h^{k-1},\partial_q^\tau e_u^k) - \sum_{i=0}^q \frac1\tau \delta_i \big(\mh^{k-i}(e_u^{k-i} , r_h^k - \eta r_h^{k-1}) \big)\Big] \nonumber \\
			&\ + \tau \sum_{k=q+1}^n \ah^k(e_w^k,r_h^k - \eta r_h^{k-1}) - \mh^k(d_1^k, r_h^k - \eta r_h^{k-1}) + c \kappa_h^n\nonumber \\
			=:&\  (I)  + (II),
		\end{align*}
		where we estimate $(II)$ as in \eqref{eq:estCIII} to
		\begin{align*}
			(II) \leq \rho \tau \sum_{k=q+1}^n | e_w^{k} |_{\ah^k}^2 + c \kappa_h^n,
		\end{align*}
		and $(I)$ is estimated by Lemma~\ref{lem:SurfDifference} and the \emph{preparations} to
		\begin{align}
			(I) =&\ \tau \sum_{k=q+1}^n \sum_{i=0}^q \delta_i \frac1\tau \Big[\mh^k(r_h^k - \eta r_h^{k-1}, e_u^{k-i}) - \mh^{k-i}(r_h^k - \eta r_h^{k-1}, e_u^{k-i})\Big] \nonumber \\
			\leq&\ c \kappa_h^n. \label{eq:Surf Helper partial difference}
		\end{align}
		In total this yields the estimate
		\begin{align}
			\tau \sum_{k=q+1}^n R_3^k \leq&\ \rho \tau \sum_{k=q+1}^n | e_w^{k} |_{\ah^k}^2 + c \kappa_h^n. \label{eq:Surf C R3}
		\end{align}
		
		The term $R_4$ is estimated similar to \eqref{eq:estCIV} by a discrete product rule, by using the following identity:
		\begin{align*}
			\mh^k(d_2^k,\partial_q^\tau e_u^k) =&\ \sum_{j=0}^q \sigma_j \mh^k (d_2^k,\partial^\tau e_u^{k-j}) \nonumber \\
			=&\ \sum_{j=0}^q \sigma_j \frac1\tau \Big(\mh^k(d_2^k,e_u^{k-j}) - \mh^{k-1}(d_2^{k-1},e_u^{k-1-j})\Big) \nonumber \\
			&\ - \sum_{j=0}^q \sigma_j \frac1\tau \Big(\mh^k(d_2^k,e_u^{k-1-j}) - \mh^{k-1}(d_2^{k},e_u^{k-1-j})\Big) \nonumber \\
			&\ - \sum_{j=0}^q \sigma_j \mh^{k-1} (\partial^\tau d_2^k, e_u^{k-i-1}).
		\end{align*}
		This leads to the exact same terms with the exception of the term involving the surface discrepancy, which is estimated by a variant of Lemma~\ref{lem:SurfDifference} using the dual norm instead of an $L^2$-$L^2$-estimate. Finally we obtain the estimate
		\begin{align}
			\tau \sum_{k=q+1}^n R_4^k \leq&\ \rho \sum_{k=n-q+1}^n \| e_u^{k} \|_{H^1(\Ga_h^k)}^2 + c \kappa_h^n. \label{eq:Surf C R4}
		\end{align}
		
		The terms in $R_5$ and $R_6$ are again estimated by Lemma~\ref{lem:SurfDifference}, following the lines of \eqref{eq:Surf Helper derivative} and \eqref{eq:Surf Helper partial difference} and yield
		\begin{align}
			\tau \sum_{k=q+1}^n R_5^k \leq&\  c \kappa_h^n, \label{eq:Surf C R5} \\
			\tau \sum_{k=q+1}^n R_6^k \leq&\ \rho \tau \sum_{k=q+1}^n \| e_w^{k} \|_{H^1(\Ga_h^k)}^2 + c \kappa_h^n. \label{eq:Surf C R6}
		\end{align}
		
		Inserting the estimates \eqref{eq:Surf C L1}, \eqref{eq:Surf C R1}, \eqref{eq:Surf C R2}, \eqref{eq:Surf C R3}, \eqref{eq:Surf C R4}, \eqref{eq:Surf C R5} and \eqref{eq:Surf C R6} into \eqref{eq:Surf Eq tested}, using the $G$-norm estimates \eqref{eq:GnormEst} and absorbing terms we arrive at
		\begin{align}
			\sum_{k=n-q+1}^n |\nabla_{\Ga_h} e_u^k|_{L^2(\Ga_h^k)}^2 + \tau \sum_{k=q+1}^n | \nabla_{\Ga_h} e_w^k|_{L^2(\Ga_h^k)}^2 \leq&\  \rho \sum_{k=n-q+1}^n \|e_u^k\|_{H^1(\Ga_h^k)}^2 \nonumber \\
			&\ + \rho \tau \sum_{k=q+1}^n \| e_w^k\|_{H^1(\Ga_h^k)}^2 + c \kappa_h^n. \label{eq:Surf ende Part C}
		\end{align}
		
		\textbf{Part D}: As in \textbf{Part D} of Proposition~\ref{prop:stability}, we test the equations \eqref{eq:SurfError} by $1$ to conclude an \emph{almost mass conservation} of the error terms and then use the Poincar\'e--Wirtinger inequality on the surface \cite[Theorem~2.12]{DziukElliott2013}.
		
		For the first equation \eqref{eq:SurfError1}, testing with $1 \in S_h^n$ yields, as in \eqref{eq:MassEq1}, 
		\begin{align*}
			\partial_q^\tau \m^n = \bfd^n:= \mh^n(d_1^n,1), 
		\end{align*}
		with $\m^n:= \mh^n(e_u^n,1)$. And we therefore obtain, following exactly the lines in \textbf{Part D} of Proposition~\ref{prop:stability}, the first Poincar\'e-type estimate:
		\begin{align}
			\| e_u^n \|_{H^1(\Ga_h^n)}^2 \leq | \nabla_{\Ga_h} e_u^n |_{L^2(\Ga_h^n)}^2 + c \kappa_h^n. \label{eq:Surf Poin 1}
		\end{align}
		
		Testing \eqref{eq:SurfError2} with $1 \in S_h^n$ yields
		\begin{align*}
			\tau \sum_{k=q+1}^n | \mh^k(e_w^k,1)|^2 \leq c \kappa_h^n,
		\end{align*}
		which results in the second Poincar\'e-type estimate
		\begin{align}
			\tau \sum_{k=q+1}^n \| e_u^k \|_{H^1(\Ga_h^k)}^2 \leq \tau \sum_{k=q+1}^n | \nabla_{\Ga_h} e_u^k |_{L^2(\Ga_h^k)}^2 + c \kappa_h^n. \label{eq:Surf Poin 2}
		\end{align}
		
		\textbf{Part E}: Inserting the Poincar\'e-type estimates \eqref{eq:Surf Poin 1} and \eqref{eq:Surf Poin 2} into \eqref{eq:Surf ende Part C} and absorbing terms by choosing the parameter $\rho>0$ small enough and Gr\"onwall's inequality, we arrive at the desired stability estimate:
		\begin{align*}
			\| e_u^n \|_{H^1(\Ga_h^n)}^2 + \tau \sum_{k=q+1}^n \| e_u^k \|_{H^1(\Ga_h^k)}^2 \leq C (\calE_I^{q-1} + D_h^q).
		\end{align*} 
		
		\textit{(The special case $n=q$)} As in the proof of Proposition~\ref{prop:stability}, the initial case $n=q$ has to be obtained in a slightly different way, see \eqref{eq:qest}, since there is no equation for $n-1=q-1$ available. We obtain
		\begin{align}
			\| e_u^q \|_{H^1(\Ga_h^n)}^2 + \tau \| e_w^q \|_{H^1(\Ga_h^n)}^2 \leq C (\calE_I^{q-1} + D_h^q),
		\end{align}
		by an analogous argument.
		
		\textit{(Proving $t_{\max}=T$)} As noted in the \emph{preparations} it is left to show that $t_{\max} =T$. This can be obtained by the exact same proof as in Proposition~\ref{prop:stability}.
	\end{proof}
	
	\subsection{Consistency estimates}
	
	\begin{prop} \label{prop:Surf consistency}
		Let $(u,w)$ solve \eqref{eq:CHevolSurf} and satisfy the regularity assumptions \eqref{eq:SurfregularityAss}. Then we have, for $q \leq k \leq n$ and $1\leq i \leq q-1$,
		\begin{alignat*}{3}
			\| d^k_1 \|_{\ast,k} + \| d^k_2 \|_{\ast,k} + \| \partial^\tau_q d^k_2 \|_{\ast,k} \leq&\  C (h^2 + \tau^q), \\
			\| \partial^\tau d^i_2 \|_{\ast,i} \leq&\   C (h^2 + \tau^q) .
		\end{alignat*}
		With these estimates we obtain the bound on the full defect term, defined in \eqref{eq:SurfDefDhn},
		\begin{align*}
			D_h^n \leq C (h^2 + \tau^{q})^2.
		\end{align*}
	\end{prop}
	\begin{proof}
		The defect terms are established similar to Section~\ref{section:Consistency} on the consistency analysis, with the geometric, Ritz map and temporal error estimates, combining the techniques developed in \cite[Section~7]{BullerjahnKovacs2024} with semi-discrete defect estimates in \cite[Section~6]{BeschleKovacs2022} and the adaptation to time dependent bilinear forms as in \cite[Section~11]{KovacsLiLubich2019} and similar to the fully discrete defect estimates in \cite{ElliottSales2024}. We split the defect into a temporal part and the spatial defect, where for the temporal part we use Lemma~\ref{Lemma:PeanoKernel} together with Lemma~\ref{lem:SurfDifference} and \cite[Proposition~2.13]{ElliottSales2024} to obtain a bound of optimal order in time. Together with the spatial defect bounds, see \cite[Section~6]{BeschleKovacs2022}, this yields the final bound, see also \cite{ElliottSales2024}. The estimates for the discrete derivative of the defect term $d_2$ are again obtained using the techniques developed in \cite[Section~7]{BullerjahnKovacs2024}, using Lemma~\ref{lem:SurfDifference} and \cite[Proposition~2.13]{ElliottSales2024} to deal with the surface discrepancy, and the spatial defects in \cite[Section~6]{BeschleKovacs2022}.
	\end{proof}

	\subsection{Optimal-order convergence result}
	
	We combine the stability result, Proposition~\ref{prop:Surfstability}, and the consistency results, Proposition~\ref{prop:Surf consistency}, as in Section \ref{section:MainProof}, to get the following optimal-order convergence result.
	
	%

	\begin{thm}\label{thm:SurfOptimalOrderErrorEst}
		Let $1\leq q \leq 5$ and $(u,w)$ be a sufficiently smooth solution (see \eqref{eq:SurfregularityAss}) of the Cahn--Hilliard equation on an evolving surface \eqref{eq:CHevolSurf}, with the nonlinear potential satisfying \eqref{eq:SurfRegAssPot}. Let $\Ga(t)$ be a $C^3$-surface, then there exists $h_0>0$ such that for all $h\leq h_0$ and $\tau >0$, satisfying the mild CFL condition $\tau^q\leq C_0 h^2$ (where $C_0>0$ can be chosen arbitrarily), the error between the solution $(u,w)$ and the fully discrete solution $(u_h^n,w_h^n)$ of \eqref{eq:SurfNumScheme}, using linear bulk--surface finite element discretizations and linearly implicit backward difference time discretization of order $q=1,\dots,5$, satisfies the optimal-order error estimates, for $q\tau \leq t_n:=n\tau \leq T$,
		\begin{align*}
			\| (u^n_h)^\ell - u(t_n) \|_{L^2(\Ga(t_n))} + h \| (u^n_h)^\ell - u(t_n) \|_{H^1(\Ga(t_n))} \leq&\ C (h^2+\tau^q),  \\
			\bigg( \tau \sum_{k=q}^n \| (w^k_h)^\ell - w(t_k) \|_{L^2(\Ga(t_n))}^2 + h \| (w^k_h)^\ell - w(t_k) \|_{H^1(\Ga(t_n))}^2 \bigg)^{\frac{1}{2}} \leq&\ C (h^2+\tau^q),
		\end{align*} 
		provided the starting values $u_h^i \in S_h^i$ are $O(\tau^q+h^2)$ accurate. The constant $C>0$ depends on the Sobolev norms of the exact solution and exponentially on the final time $T$, but is independent of $h$ and $\tau$.
	\end{thm}
	
	Sufficient regularity assumptions for Theorem~\ref{thm:SurfOptimalOrderErrorEst} are: the nonlinear potential and its derivative 
	\begin{equation}
		F_\Ga^{(k)} \text{ is locally Lipschitz continuous for }k=0,\dots,3, \label{eq:SurfRegAssPot}
	\end{equation}
	and the exact solutions satisfy 
	\begin{equation} \label{eq:SurfregularityAss}
		\begin{aligned}
			u,\matdev u, \dots, (\matdev)^{q+1} u \in H^2(\Ga(t)), \\
			w, \matdev w, (\matdev)^2 w  \in H^2(\Ga(t)),
		\end{aligned}
		\quad \text{uniformly in } t \in [0,T],
	\end{equation}
	with the surface velocity field satisfying $V, \matdev V, \dots, (\matdev)^q V \in W^{1,\infty}(\Ga(t))$, uniformly in $t \in [0,T]$.
	\begin{rem}
		(i) Note that the function spaces appearing in Theorem~\ref{thm:SurfOptimalOrderErrorEst} and \eqref{eq:SurfregularityAss} are \emph{evolving} Bochner spaces as formally introduced in \cite{AlphonsoElliottStinner2015}.
		
		(ii) Comparing Theorem~\ref{thm:SurfOptimalOrderErrorEst} to the main result in \cite[Theorem~4.2]{ElliottSales2024} for a fully implicit backward Euler time discretization one notices that our linearly implicit method does not suffer from the difficulties appearing in the stability, and allows a less restrictive choice of the nonlinear potential, cf. \eqref{eq:SurfRegAssPot}. While the regularity assumptions in \eqref{eq:SurfregularityAss} are similar compared to \cite[Equation~(4.1)]{ElliottSales2024}. Note that neither of the regularity assumptions are covered by the well-posedness result for the exact solution in \cite[Theorem~4.14]{CaetanoElliott2021}, even though the spatial regularity is achieved.
	\end{rem}

	\section{Numerical experiments}
	\label{section:NumericalExperiments}
	In this section we illustrate and complement our theoretical results with some numerical experiments concerning the order of convergence in space and time, as well as some experiments on the evolution of solutions for different initial data, repeating some computations from \cite{KnopfLam2020} and \cite{KnopfLamLiuMetzger2021}, allowing for a direct comparison. We present the following experiments:
	\begin{enumerate}
		\item Computation of convergence rates by constructing manufactured solutions on the unit ball, cf. \cite{BullerjahnKovacs2024}.
		\item Evolution of an elliptically shaped droplet for various values $K \in [0,\infty)$ and $L \in (0,\infty]$ extending on the numerical experiments in \cite{KnopfLamLiuMetzger2021}.
		\item Evolution of a randomly distributed initial data for an affine linear transmission condition, comparing to the experiments in \cite{KnopfLam2020}.
	\end{enumerate}
	In all examples we use quadrature rules of sufficiently high order to avoid issues of numerical integration.
	\subsection{Numerical experiments on convergence rates}
	We analyse the rate of convergence for the bulk--surface Cahn--Hilliard system with dynamic boundary conditions \eqref{eq:CHreact}, with the choice of parameters $\epsilon=\delta=\kappa=m_\Om=m_\Ga=1$ and $\alpha=1, \beta=1$, on the two dimensional unit disk $\Om$ with a manufactured solution. This means we add additional inhomogeneities on the right-hand side of each equation, to obtain as the exact solution:
	\begin{align*}
		u(t,x_1,x_2)= e^{-t}x_1x_2, \qquad \psi(t,x_1,x_2)= e^{-t}x_1x_2, \\ \mu(t,x_1,x_2)= e^{-t}x_1x_2, \qquad \theta(t,x_1,x_2)= e^{-t}x_1x_2,
	\end{align*}
	for the nonlinear double well potentials
	\begin{align}
		F_\Om(s)=&\ \frac{1}{4}(s^2-1)^2, \qquad F_\Ga(s)= \frac{1}{4}(s^2-1)^2. \label{eq:DoubleWell}
	\end{align} 
	
	\begin{rem}
		By a minimal extension of problem \eqref{eq:CHreact}, with an additional inhomogeneity for each equation, this fits into the framework of the theoretical result Theorem~\ref{thm:OptimalOrderErrorEst}. The error equations are not changed by this modification, only some simple defect terms are added, which are due to the geometric approximation, and we need $L^2(\Om)$-regularity in the bulk and $L^2(\Ga)$-regularity on the surface for the inhomogeneities, to obtain the analogous regularity result.
	\end{rem} 
	
	We use the linear bulk--surface finite elements in space and the 2-step BDF method to discretize this problem in the time interval $[0,1]$, with starting values taken from the exact solution.  
	
\begin{figure}[!t]
		\centering\includegraphics[width=\textwidth]{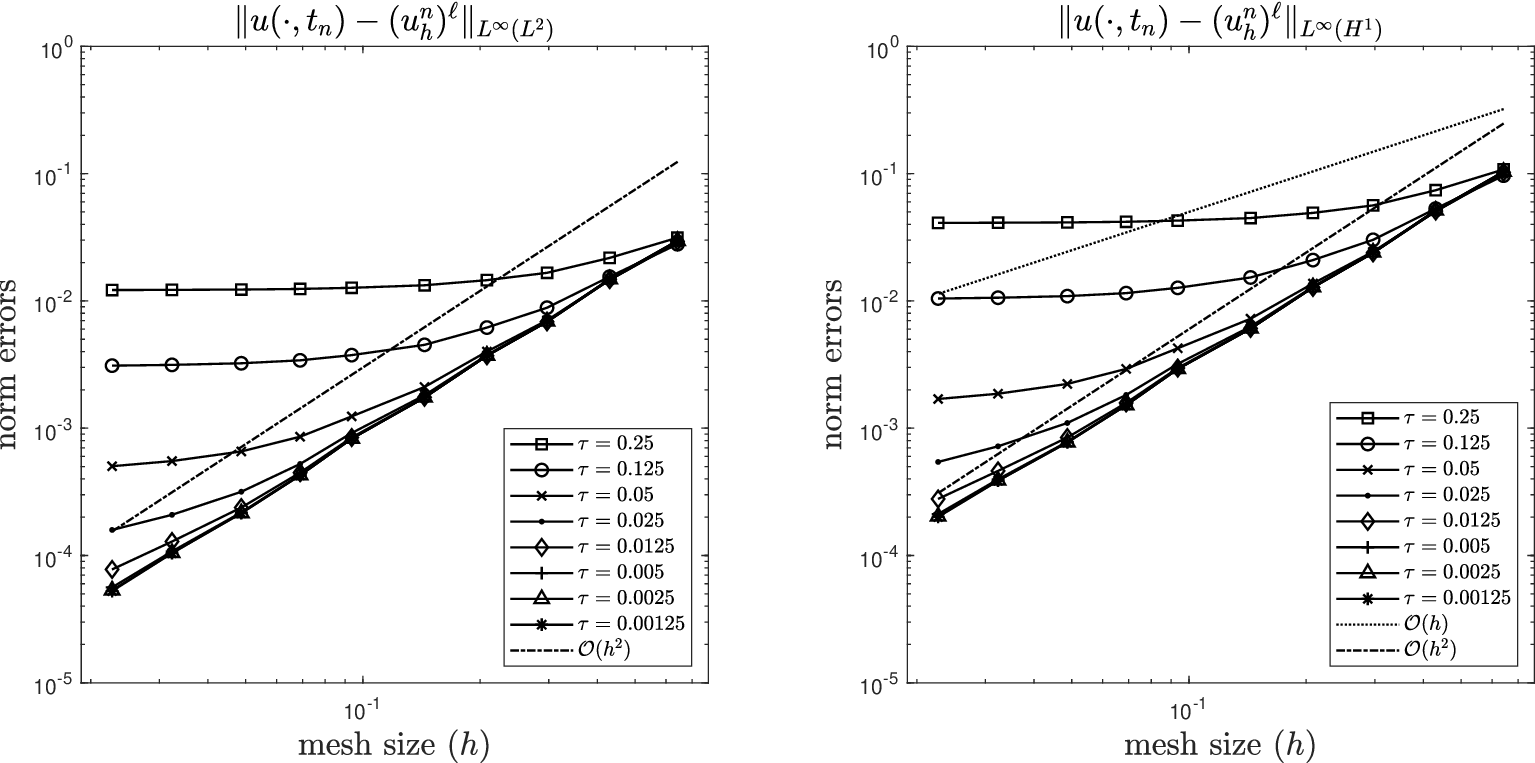}
		\caption{Spatial convergence plot for the linear bulk--surface FEM/BDF2 approximation to the bulk--surface Cahn--Hilliard system with dynamic boundary conditions for $K=0$ and $L=100$}
		\label{fig:SpatialConvergence1}
	\end{figure}

\begin{figure}[!t]
		\centering\includegraphics[width=\textwidth]{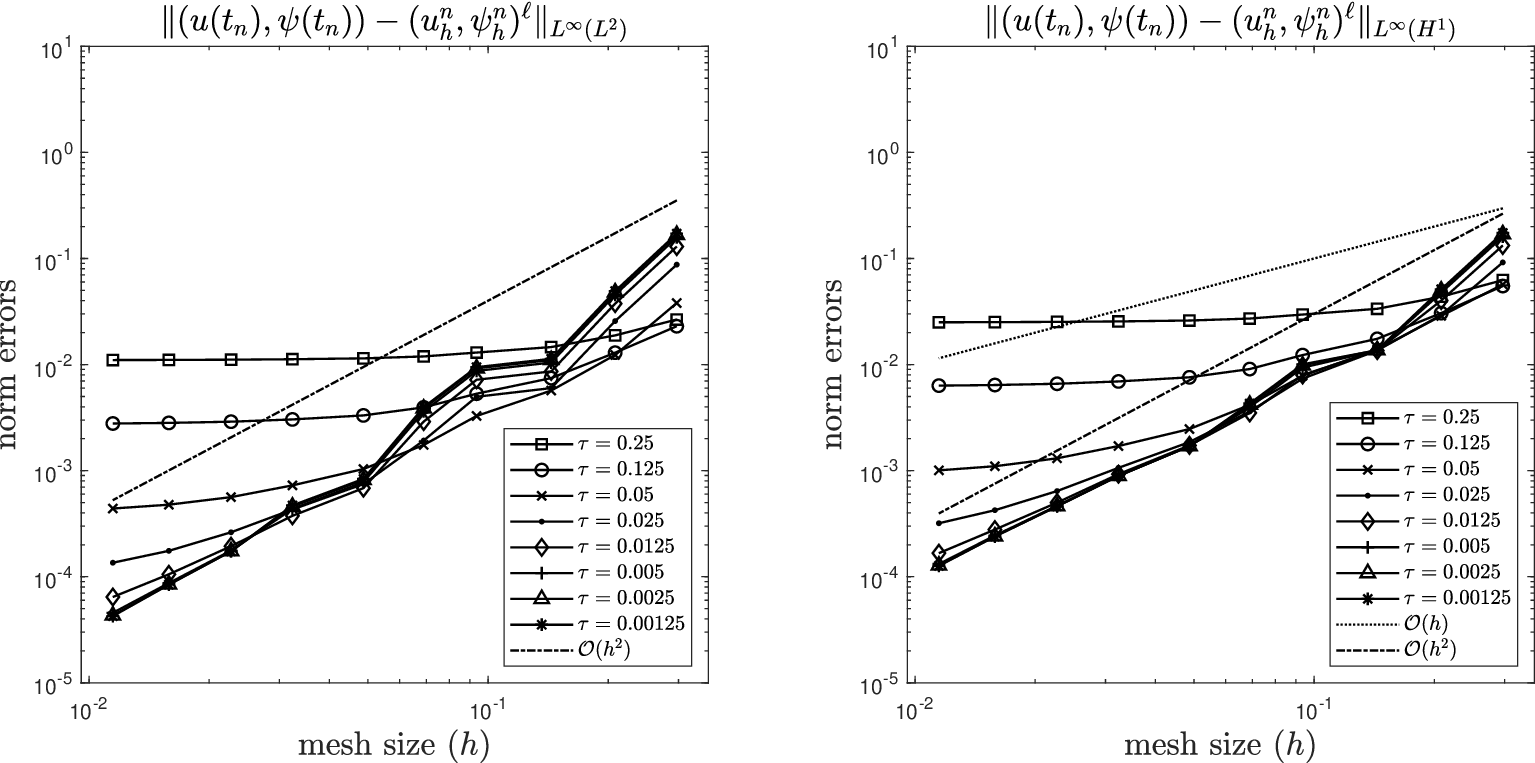}
		\caption{Spatial convergence plot for the linear bulk--surface FEM/BDF2 approximation to the bulk--surface Cahn--Hilliard system with dynamic boundary conditions for $K=10$ and $L=0.01$.}
		\label{fig:SpatialConvergence2}
	\end{figure}

	We choose two sets of extremal parameters $K=0$, $L=100$, and $K=10$, $L=0.01$, and display in Figure~\ref{fig:SpatialConvergence1}, respectively, Figure~\ref{fig:SpatialConvergence2}, the $L^2$- and $H^1$-norm errors for the phase fields $u$ and $\psi$, in a logarithmic plot for meshes with degrees of freedom $2^k \cdot 20$ for $k \in \{2,...,11\}$ and $\tau \in \{0.25, 0.125,0.05,0.025,0.0125,0.005, 0.0025,0.00125\}$. We notice that for a small enough time step size the spatial error dominates and is of order $h^2$, as expected from Theorem~\ref{thm:OptimalOrderErrorEst} for the $L^\infty([0,1],L^2)$-norm and even one order higher as the theoretical result in the case of the $L^\infty([0,1],H^1)$-norm (due to the super-convergence when comparing with the interpolation of the exact solution, cf. \cite[Section~9.1]{HarderKovacs2022}). 
	
\begin{figure}[!t]
		\centering\includegraphics[width=\textwidth]{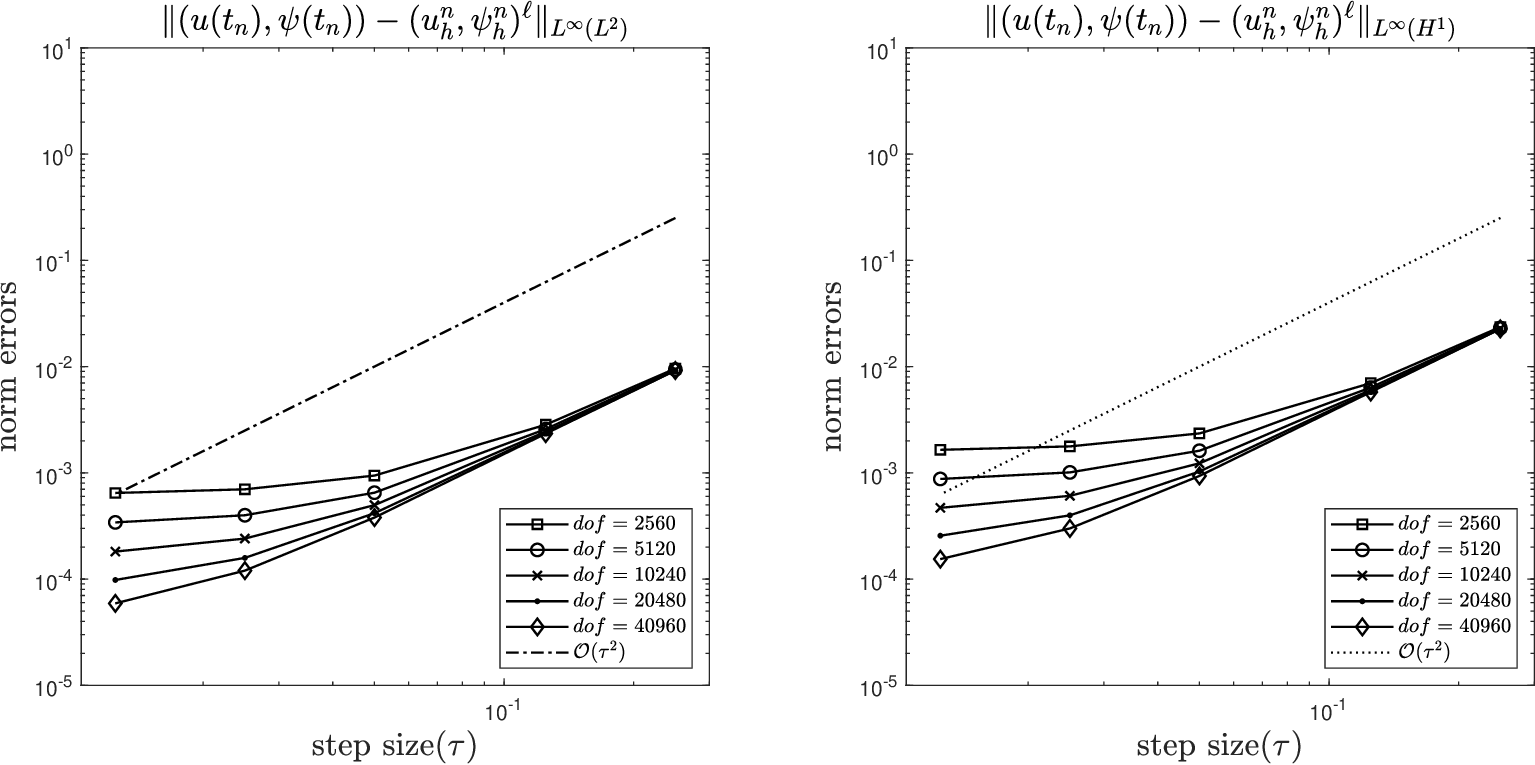}
		\caption{Temporal convergence plot for the linear bulk--surface FEM/BDF2 approximation to the bulk--surface Cahn--Hilliard system with dynamic boundary conditions for $K=0$ and $L=100$}
		\label{fig:TemporalConvergence1}
	\end{figure}

\begin{figure}[!t]
		\centering\includegraphics[width=\textwidth]{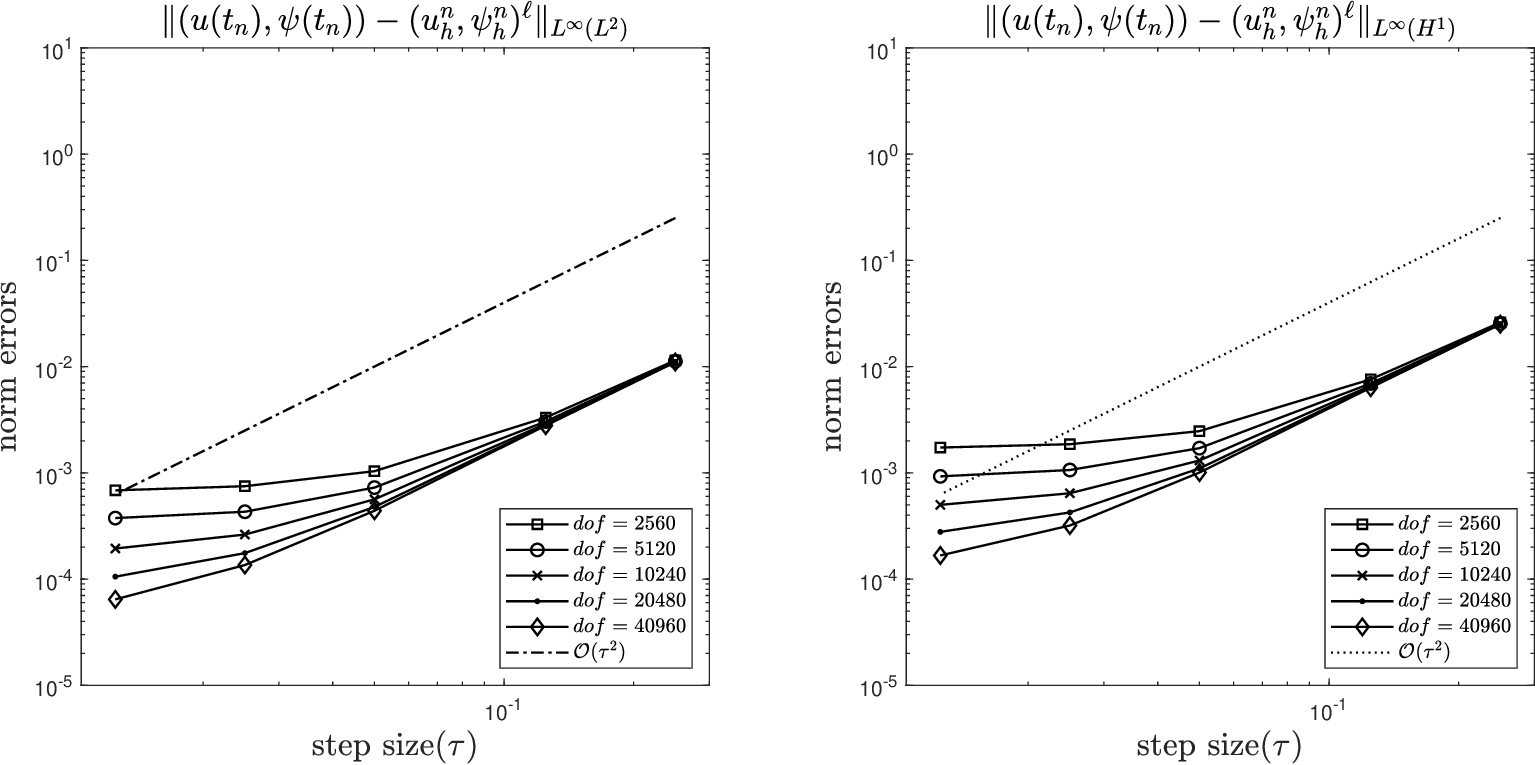}
		\caption{Temporal convergence plot for the linear bulk--surface FEM/BDF2 approximation to the bulk--surface Cahn--Hilliard system with dynamic boundary conditions for $K=10$ and $L=0.01$}
		\label{fig:TemporalConvergence2}
	\end{figure}

	In Figures \ref{fig:TemporalConvergence1} and \ref{fig:TemporalConvergence2}, the same errors are displayed for parameters $K=0$ and $L=100$, respectively, $K=10$ and $L=0.01$, now against the time step size $\tau  \in \{0.25, 0.125, 0.05, 0.025, 0.0125\}$, and with variable mesh sizes, corresponding to degrees of freedom $2^k \cdot 20$ for $k \in \{7,\dotsc,11\}$. We again observe that, for a small enough mesh size $h$, the temporal error dominates and admits the convergence rate $\tau^2$ for the $2$-step BDF method, which is the expected convergence by our theoretical result Theorem~\ref{thm:OptimalOrderErrorEst}. 
	
	For other parameter values $K,L$  we have observed identical convergence behavior.
	
	\subsection{Evolution of an elliptically shaped droplet with respect to $K$ and $L$}
	
\begin{figure}[!t]
		\centering\includegraphics[width=0.5\textwidth]{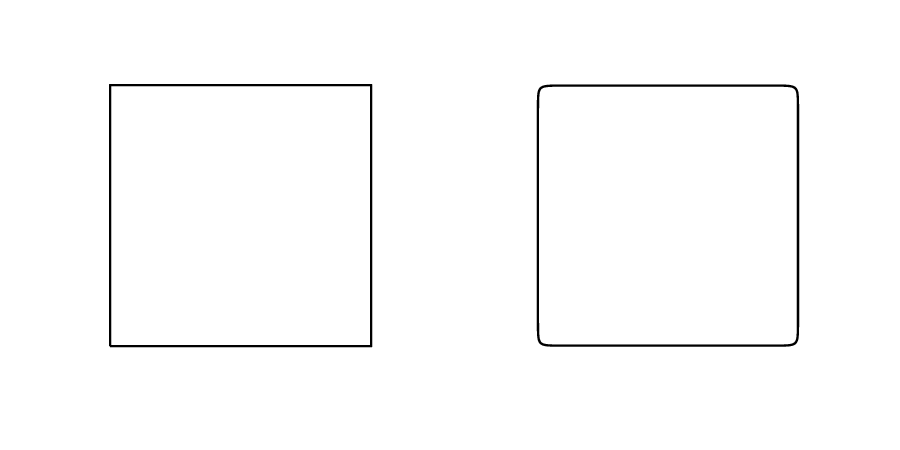}
		\caption{The unit square (left) and the approximated rounded square (right).}
		\label{fig:RoundSquare}
	\end{figure}

	We consider the bulk--surface Cahn--Hilliard system with dynamic boundary conditions \eqref{eq:CHreact}, and investigate the behavior of the solution $u$ for different transmission rates $K \in [0,\infty)$ and reaction rates $L \in (0,\infty]$ for an elliptically shaped droplet. In order to compare to the experiments conducted in \cite[Section~6]{KnopfLamLiuMetzger2021} on the unit square $(0,1)^2$, yet still satisfying the assumptions of our theoretical results, we choose the domain $\Om_\omega$ as a approximate rounded square with a $C^\infty$-boundary. The rounded square $\Om_\omega$ is given by the level set function
	\begin{equation*}
		\phi_\omega (x,y) = \omega \log(e^{\frac{x-1}{\omega}} + e^{\frac{-x}{\omega}} + e^{\frac{y-1}{\omega}} + e^{\frac{-y}{\omega}}),
	\end{equation*}
	where $\omega=0.015$ in our examples, which leads to a maximal distance of approximately $\omega \log(2)$ for the two boundaries, cf. Figure~\ref{fig:RoundSquare}.
	
	We set the parameters in the equation to $\epsilon=\delta=0.02$ and $\beta=\alpha=m_\Om=m_\Ga=\kappa=1$, and define the non-linear potentials as the double well potentials
	
	\begin{align*}
		F_\Om(s)=&\ \frac{1}{8}(s^2-1)^2, \qquad F_\Ga(s)= \frac{1}{8}(s^2-1)^2.
	\end{align*} 
	
	We approximate the solution by our numerical scheme, combining the linear finite element method and the linearly implicit BDF method of order $2$, with a triangular mesh consisting of $20500$ nodes and $\tau=10^{-5}$. As the initial data $(\bfu^0,\bfpsi^0)$ we choose an elliptical shaped droplet with barycenter at $(0.1,0.5)$, a maximal horizontal elongation of $0.6814$, and a maximal vertical elongation of $0.367$, see also \cite[Section~6]{KnopfLamLiuMetzger2021}. The second required initial data $(\bfu^1,\bfpsi^1)$ is computed by a linearly implicit BDF step of order $1$.

	\begin{figure}[!t]
		\centering\includegraphics[width=\textwidth]{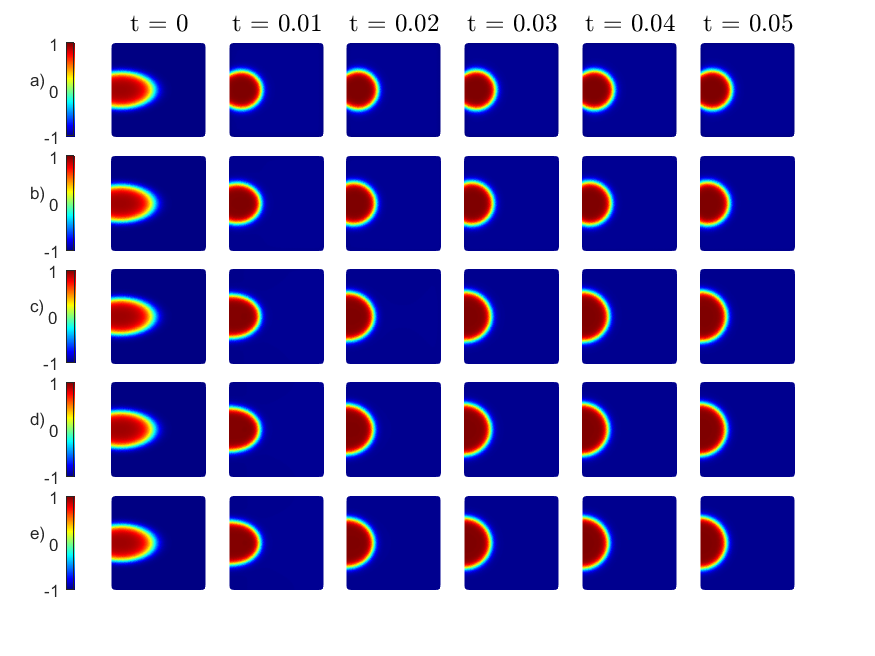}
		\caption{Evolution plot for the linear bulk--surface FEM/BDF2 approximation to \eqref{eq:CHreact}, and time steps $t=0,0.01,0.02,0.03,0.04,0.05$, for the parameter $a) ~ K=10^{-5}$, $b) ~ K=0.1$, $c) ~ K=1$, $d) ~ K=10$ and $e) ~ K=10^5$.}
		\label{fig:EvolDroplet1}
	\end{figure}
	
	In Figure~\ref{fig:EvolDroplet1} we fix $L=\infty$ and observe that over time the elliptical droplet tries to attain a circular shape with constant mean curvature, but with the restriction that the contact area with the boundary only changes depending on the parameter $K \in \{10^{-5}, 0.1,1,10,10^5\}$: The bigger the value of $K$, the more the contact area increases. Additionally as $K \to \infty$ we observe that the contact angle of the interfacial region and the boundary approaches $\pi/2$ as expected from the Neumann boundary conditions.
	
\begin{figure}[!t]
		\centering\includegraphics[width=\textwidth]{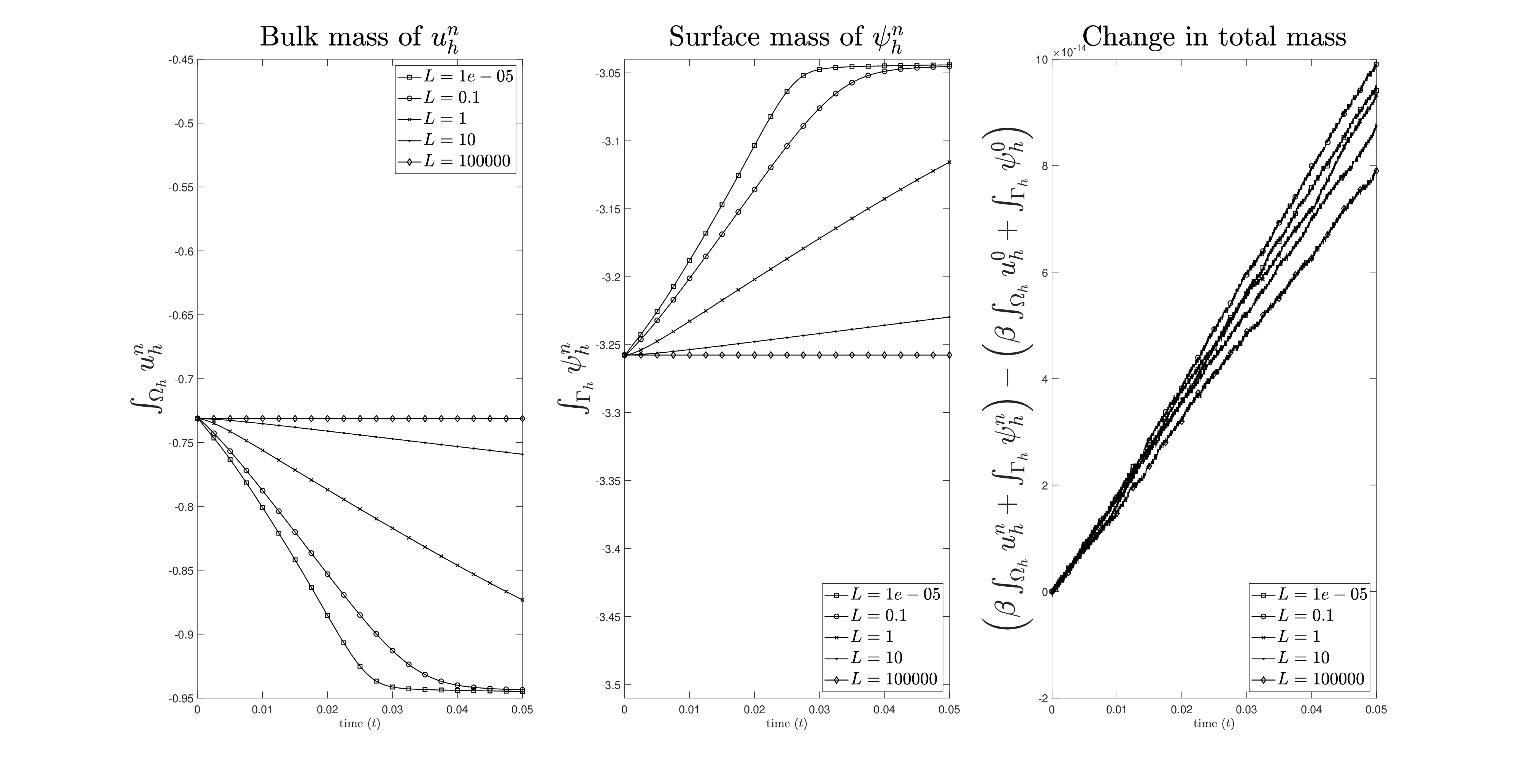}
		\caption{Evolution of the bulk, surface and total mass for the linear bulk--surface FEM/BDF2 approximation to the Cahn--Hilliard equation with reaction rate dependent dynamic boundary conditions for an elliptically shaped droplet (note axis scaling).}
		\label{fig:MassDroplet}
	\end{figure}

\begin{figure}[!t]
		\centering\includegraphics[width=\textwidth]{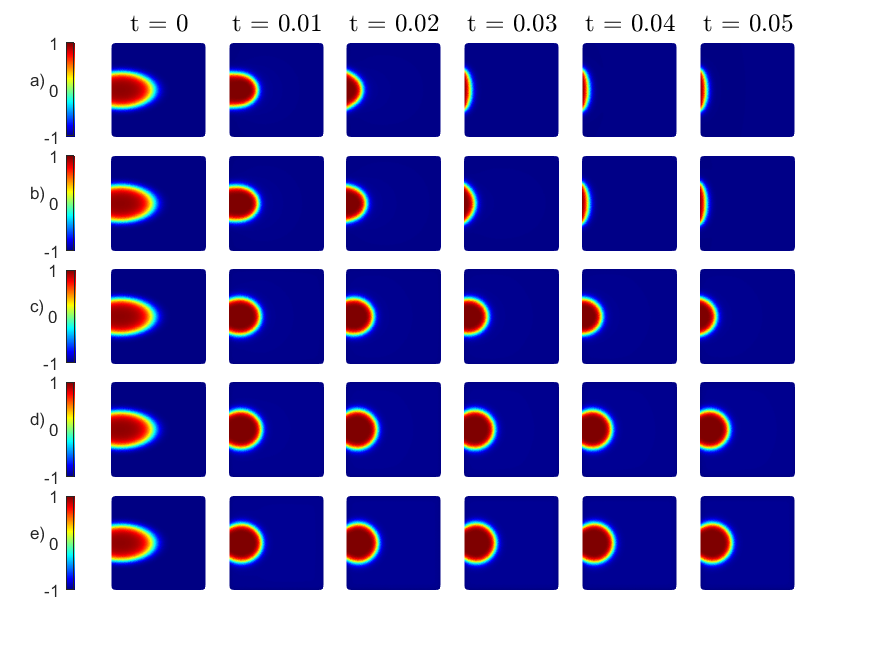}
		\caption{Evolution plot for the linear bulk--surface FEM/BDF2 approximation to \eqref{eq:CHreact}, and time steps $t=0,0.01,0.02,0.03,0.04,0.05$, for the parameter $a) ~ L=10^{-5}$, $b) ~ L=0.1$, $c) ~ L=1$, $d) ~ L=10$ and $e) ~ L=10^5$.}
		\label{fig:EvolDroplet2}
	\end{figure}

	When fixing $K=0$, Figure \ref{fig:MassDroplet} shows the evolution of the bulk, surface and total mass for the droplet. As we know, the exact solution satisfies a conservation of mass for the bulk and surface separately for $L=\infty$ and a conservation of the total mass $\beta \int_\Om u + \int_\Ga 1$ for all $L \in [0,\infty)$, see \eqref{eq:MassConservation}, and we observe, as expected, that our numerical solution replicates this exactly. Over time the elliptical droplet tries to attain a circular shape with constant mean curvature, cf. Figure \ref{fig:EvolDroplet2}, but with the restriction that the contact area only changes depending on the conservation of surface mass, determined by $L \in \{10^{-5}, 0.1,1,10,10^5\}$: The smaller the value of $L$, the more the contact area increases and bulk mass decreases, comparing to the experiments in \cite{KnopfLamLiuMetzger2021}.
	
	\subsection{Evolution of a randomly distributed initial data for an affine linear transmission condition with $K=10^{-5}$ and $L=\infty$}
	
	We consider the bulk--surface Cahn--Hilliard system with dynamic boundary conditions \eqref{eq:CHreact} on the rounded square $\Om_\omega$ for parameters $K=10^{-5}$ and $L=\infty$, with an \emph{affine linear transmission condition} $$\epsilon K \partial_\nu = \alpha_1 \psi + \alpha_2 -u$$ in place of \eqref{eq:robintypecond1}, cf. \cite{KnopfLam2020}.
	
	\begin{rem}
		By a minimal extension the problem \eqref{eq:CHreact}, with an affine linear transmission condition, fits into the framework of the theoretical result Theorem~\ref{thm:OptimalOrderErrorEst}. The error equations are not changed by this modification, only a simple defect term is added, which is due to the geometric approximation, and can be treated easily with the established methods, to obtain the analogous regularity result.
	\end{rem}
	
	We would like to investigate the behavior of the solution for different values $\alpha_2 \in \{-0.3,0.3,0.9,1.2 \}$, comparing, by the choice of $K=10^{-5}$, to the findings in \cite[Section~7.2.1]{KnopfLam2020}. We set the parameters in the equation to $\epsilon=\delta=0.02$ and $\alpha_1=m_\Om=m_\Ga=\kappa=1$, and define the non-linear potentials as the double well potentials, see \eqref{eq:DoubleWell}.
	
	We approximate the solution by our numerical scheme, combining the linear finite element method and the linearly implicit BDF method of order $2$, with a triangular mesh consisting of $20500$ nodes and $\tau=10^{-6}$. 
	
	As the initial data $\bfu^0$ we choose a uniformly distributed random value in $[0.3,0.5]$ (this was generated by the \textsc{Matlab}-function 'rand' with the Mersenne Twister generator using seed 42), such that $\langle u_0 \rangle_\Om=0.3997$ and $\langle u_0 \rangle_\Ga=0.4004$, and on the boundary we let $\bfpsi^0=\alpha_1^{-1}(\bfu^0|_{\Ga_h} - \alpha_2)$, in accordance with the Dirichlet boundary condition for $K=0$. The second required initial data $(\bfu^1,\bfpsi^1)$ is computed by a linearly implicit BDF step of order $1$ with smaller time step size. 
	
\begin{figure}[!t]
		\centering\includegraphics[width=\textwidth]{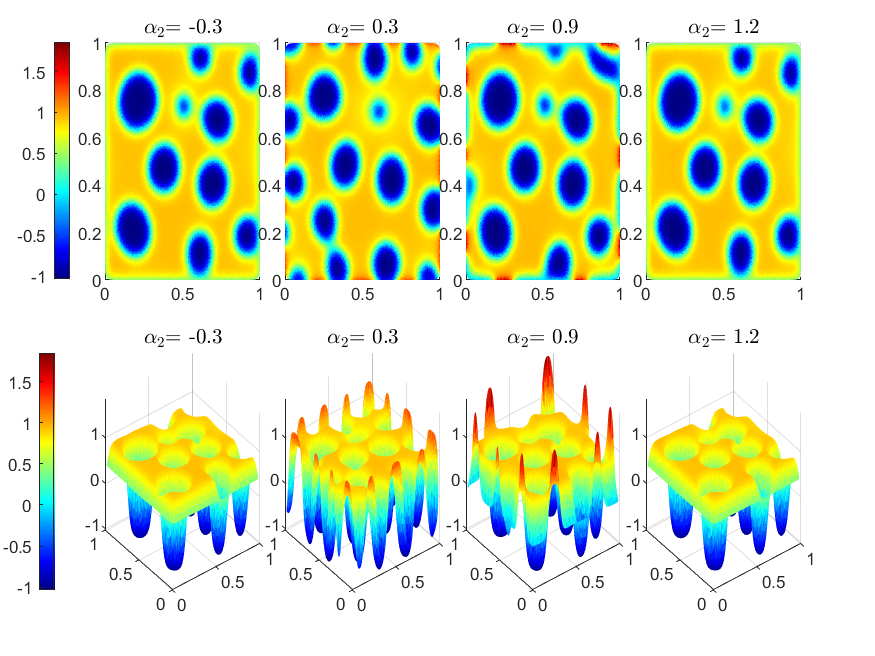}
		\caption{The linear bulk--surface FEM/BDF2 approximation $\bfu^n$ for \eqref{eq:CHreact}, with the parameters $K=10^{-5}$, $L=\infty$, $\alpha_1=1$ and $\alpha_2 \in \{-0.3,0.3,0.9,1.2\}$, at iteration $n=2001$.}
		\label{fig:Evolbeta}
	\end{figure}

	We observe in Figure~\ref{fig:Evolbeta} the phenomenon that the behavior of the bulk phase field at the boundary depends heavily on the parameter $\alpha_2$. In the cases $\alpha_2=0.3,0.9$ we observe the Cahn--Hilliard phase separation dynamics on the boundary, in contrast to the cases $\alpha_2=-0.3,1.2$ where we observe near constant trace values, close to $\langle u_0 \rangle_\Ga$. These results replicate the experiments done in \cite[Section~7.2.1]{KnopfLam2020} for a finite element/backward Euler scheme, as indicated therein, we do not have a robust method to predict whether phase separation or near constant trace values will occur.
	
	\section*{Acknowledgments}
	I would like to thank Bal\'azs Kov\'acs for carefully reading the manuscript and providing many valuable comments and suggestions.

	\bibliography{C-Heq_system}

\end{document}